\input amstex\documentstyle{amsppt}  
\pagewidth{12.5cm}\pageheight{19cm}\magnification\magstep1
\topmatter
\title Unipotent representations as a categorical centre\endtitle
\author G. Lusztig\endauthor
\address{Department of Mathematics, M.I.T., Cambridge, MA 02139}\endaddress
\thanks{Supported in part by National Science Foundation grant 1303060.}\endthanks   
\endtopmatter   
\document
\define\rk{\text{\rm rk}}

\define\Irr{\text{\rm Irr}}

\define\Bpq{\Bumpeq}

\define\dL{\dot L}

\define\bco{\bar{\co}}

\define\unz{\un{\z}}
\define\unc{\un{\c}}
\define\unb{\un{\bul}}
\define\unst{\un{*}}

\define\mpb{\medpagebreak}

\define\si{\sim}

\define\sqc{\sqcup}

\define\qua{\quad}

\define\bZ{\bar Z}

\define\op{\oplus}
   
\define\part{\partial}
\define\emp{\emptyset}

\define\iy{\infty}
\define\m{\mapsto}
\define\do{\dots}

\define\lra{\leftrightarrow}

\define\sub{\subset}    
\define\bxt{\boxtimes}
\define\T{\times}
\define\ti{\tilde}
\define\nl{\newline}
\redefine\i{^{-1}}

\define\un{\underline}
\define\ov{\overline}
\define\ot{\otimes}
\define\bbq{\bar{\QQ}_l}

\define\Ad{\text{\rm Ad}}
\define\Hom{\text{\rm Hom}}

\define\tr{\text{\rm tr}}

\redefine\spa{\spadesuit}

\define\a{\alpha}

\redefine\c{\chi}
\define\g{\gamma}

\define\e{\epsilon}

\define\io{\iota}

\define\p{\pi}
\define\ph{\phi}
\define\ps{\psi}
\define\r{\rho}
\define\s{\sigma}
\redefine\t{\tau}

\redefine\l{\lambda}
\define\z{\zeta}
\define\x{\xi}

\define\vt{\vartheta}

\redefine\D{\Delta}

\redefine\L{\Lambda}
\define\Ph{\Phi}

\redefine\aa{\bold a}

\define\boc{\bold c}

\define\ee{\bold e}

\define\kk{\bold k}

\define\pp{\bold p}

\define\ww{\bold w}

\redefine\AA{\bold A}

\define\CC{\bold C}
\define\DD{\bold D}

\define\FF{\bold F}

\define\HH{\bold H}

\define\JJ{\bold J}

\define\LL{\bold L}

\define\NN{\bold N}

\define\QQ{\bold Q}

\define\ZZ{\bold Z}
\define\XX{\bold X}

\define\ca{\Cal A}
\define\cb{\Cal B}
\define\cc{\Cal C}
\define\cd{\Cal D}

\define\ch{\Cal H}

\define\cl{\Cal L}
\define\cm{\Cal M}

\define\co{\Cal O}

\define\cz{\Cal Z}
\define\cx{\Cal X}
\define\cy{\Cal Y}

\define\fA{\frak A}

\define\fD{\frak D}
\define\fE{\frak E}

\define\fL{\frak L}

\define\fS{\frak S}
\define\fT{\frak T}

\define\te{\ti e}

\define\tih{\ti h}

\define\ty{\ti y}

\define\tF{\ti F}

\define\tL{\ti L}

\define\tP{\ti P}

\define\sh{\sharp}

\define\bul{\bullet}

\define\cir{\bul}
\define\prq{\preceq}

\define\Rep{\text{\rm Rep}}

\define\BFO{BFO}
\define\DL{DL}
\define\EW{EW}
\define\ENO{ENO}
\define\ORA{L1}
\define\CSII{L2}
\define\CSIII{L3}
\define\CSIV{L4}
\define\COX{ L5}
\define\CDGVII{L6}
\define\CDGIX{L7}
\define\CDGX{L8}
\define\CEL{L9}
\define\HEC{L10}
\define\RES{L11}
\define\TC{L12}
\define\MU{Mu}
\define\SO{So}
\head Introduction\endhead
\subhead 0.1\endsubhead
Let $\kk$ be an algebraic closure of the finite field $\FF_p$ with $p$ elements. For any power $q$ of $p$ 
let $\FF_q$ be the subfield of $\kk$ with $q$ elements. Let $G$ be a reductive connected group over $\kk$,
assumed to be adjoint. Let $\cb$ be the variety of Borel subgroups of $G$. 

Let $W$ be the Weyl group of $G$ and let $\boc$ be a two-sided cell of $W$. 
Let $s\in\ZZ_{>0}$ and let $F:G@>>>G$ be the Frobenius map for an $\FF_{p^s}$-rational structure
on $G$. Let $G^F=\{g\in G;F(g)=g\}$, a finite group. Let $\Rep^\spa(G^F)$  (resp. $\Rep^\boc(G^F)$) be the 
category of representations of $G^F$ over $\bbq$ which are finite direct sums of unipotent representations 
in the sense of \cite{\DL} (resp. of unipotent representations whose associated two-sided cell (see 1.3) is
$\boc$); here $l$ is a fixed prime number invertible in $\kk$.

In the rest of this subsection we assume for simplicity that the $\FF_{p^s}$-rational structure on $G$ is 
split. The simple objects of $\Rep^\boc(G^F)$ were classified in \cite{\ORA}. The classification turns out
to be the same as that \cite{\CSIV} 
of unipotent character sheaves on $G$ whose associated two-sided cell is $\boc$. The fact that 

(a) {\it these two classification problems have the same solution}
\nl
has not until now been adequately explained.

In \cite{\TC} we have shown that the category of perverse sheaves on $G$ which are direct sums
of unipotent character sheaves with associated two-sided cell $\boc$ is naturally equivalent to the centre of
a certain monoidal category $\cc^\boc\cb^2$ of sheaves on $\cb^2$ introduced in \cite{\CEL}
for which the induced ring structure on the Grothendieck group is the $J$-ring attached to 
$\boc$, see \cite{\HEC, 18.3}. (The analogous statement for $D$-modules on a reductive group over $\CC$ was 
proved earlier in a quite different way in \cite{\BFO}.) In this paper we show that
$\Rep^\boc(G^F)$ is also naturally equivalent to the centre of $\cc^\boc\cb^2$ (see 6.3).
This implies in particular that the simple objects of $\Rep^\boc(G^F)$ are naturally in bijection with 
the unipotent character sheaves with associated two-sided cell $\boc$, which explains (a). It also implies 
that the set of simple objects $\Rep_\boc(G^F)$ is ``independent'' of the choice of $s$; in fact, as we show
in 7.1, it is also independent of the characteristic of $\kk$. It follows that to classify the unipotent
representations of $G^F$ it is enough to classify the unipotent character sheaves on $G$ in sufficiently
large characteristic; for the latter classification one can use the scheme of \cite{\RES} which uses the 
unipotent support of a character sheaf.

The methods of this paper are extensions of those of \cite{\TC}. We replace $\Rep^\boc{G^F}$ by an equivalent
category consisting of certain $G$-equivariant perverse sheaves on $G_s$, the set of all Frobenius maps
$G@>>>G$ corresponding to split $F_{p^s}$-rational structures on $G$; we view $G_s$ as an algebraic variety 
in a natural way. We construct functors $\unc_s,\unz_s$ between this category and the category
$\cc^\boc\cb^2$ which are $q$-analogues of the truncated induction and truncated restriction $\unc,\unz$ of 
\cite{\TC} and we show that most properties of $\unc,\unz$ are preserved.
We also define a truncated convolution product from our sheaves on $G_s$ and on $G_{s'}$ to our sheaves
on $G_{s+s'}$ which is analogous to the truncated convolution of character sheaves in \cite{\TC}; we also 
give a meaning for this even when $s,s'$ are arbitrary integers. The main
application of this truncated convolution product is in the case where $s'=-s$, the result of the product 
being a direct sum of character sheaves on $G$; this is used in the proof of a weak form of an adjunction 
formula between $\unc_s,\unz_s$ which is then used to prove the main result (Theorem 6.3).

\subhead 0.2\endsubhead
In this paper we also prove extensions of the results in 0.1 to the case where $F:G@>>>G$ is the
Frobenius map of a nonsplit $F_{p^s}$-rational structure. In this case the role of unipotent character
sheaves on $G$ is taken by the unipotent character sheaves on a connected component of the group of
automorphism group of $G$. Moreover, in this case the centre of $\cc^\boc\cb^2$ is replaced by a slight
generalization of the centre (the $\e$-centre) which depends on the connected component above.

Many arguments in this paper are very similar to arguments in \cite{\TC} and are often replaced by
references to the corresponding arguments in \cite{\TC}.

Our results can be extended to non-unipotent representations and non-unipotent character sheaves; this will 
be discussed elsewhere.

\subhead 0.3\endsubhead
{\it Notation.}
We assume that we are given a split $\FF_p$-rational structure on $G$ with Frobenius map $F_0:G@>>>G$. Let 
$\nu=\dim\cb$, $\D=\dim(G)$, $\r=\rk(G)$. We shall view $W$ as an indexing set for the orbits of $G$ acting 
on $\cb^2:=\cb\T\cb$ by simultaneous conjugation; let $\co_w$ be the orbit corresponding to $w\in W$ and let
$\bco_w$ be the closure of $\co_w$ in $\cb^2$. For $w\in W$ we set $|w|=\dim\co_w-\nu$ (the length of $w$). 
Let $w_{max}$ be the unique element of $W$ such that $|w_{max}|=\nu$.

As in \cite{\ORA, 3.1}, we say that an automorphism $\e:W@>>>W$ is {\it ordinary} if it leaves stable the 
set $\{s\in W;|s|=1\}$ and for any two elements $s\ne s'$ in that set which are in the same orbit of $\e$, 
the product $ss'$ has order $\le3$. Let $\fA$ be the group of ordinary automorphisms of $W$.

For $B\in\cb$, let $U_B$ be the unipotent radical of $B$. Then $B/U_B$ is independent of $B$; it is ``the'' 
maximal torus $T$ of $G$. Let $\cx$ be the group of characters of $T$.

Let $\Rep W$ be the category of finite dimensional representations of $W$ over $\QQ$; let $\Irr W$ be a set 
of representatives for the isomorphism classes of irreducible objects of $\Rep W$. 

The notation $\cd(X),\cm(X),\cd_m(X),\cm_m(X)$ is as in \cite{\TC, 0.2}. 
(When $X$ is $G,\cb,\co_w$ or $\bco_w$, the subscript ${}_m$ refers to the $\FF_{p^{s_0}}$-structure 
defined by $F_0^{s_0}$ for a sufficiently large $s_0>0$.) For $K\in\cd(X)$, $\ch^iK$,
$\ch^i_xK$, $K^i$, $K[[m]]=K[n](n/2)]$, $\fD(K)$ are as in \cite{\TC, 0.2}. For $K\in\cm_m(X)$, $gr_jK$ is 
as in \cite{\TC, 0.2}. For $K\in\cd_m(X)$, $K^{\{i\}}=gr_i(K^i)(i/2)$, is as in \cite{\TC, 0.2}.

If $K\in\cm(X)$ and $A$ is a simple object of $\cm(X)$ we denote by $(A:K)$ the multiplicity of $A$ in a
Jordan-H\"older series of $K$. The notation $C\Bpq\{C_i;i\in I\}$ is as in \cite{\TC, 0.2}.

If $X,X'$ are algebraic varieties over $\kk$, we say that a map of sets $f:X@>>>X'$ is a quasi-morphism if
for some $\FF_q$-rational structure on $X$ and $X'$ with Frobenius maps $F$ and $F'$ and some integer 
$t\ge0$, $fF^t:X@>>>X'$ is a morphism equal to $F'{}^tf$.
If, in addition, $fF=F't$ then we have well defined functors $f_!:\cd_m(X)@>>>\cd_m(X')$, 
$f^*:\cd_m(X')@>>>\cd_m(X)$ such that $f_!$ is the composition of usual functors
$(fF^t)_!(F^t)^*=(F'{}^t)^*(F'{}^tf)_!$ and $f^*$ is the composition of usual functors
$(F^t)_!(fF^t)^*=(F'{}^tf)^*(F'{}^t)_!$.
The usual properties of $f_!,f^*$ for morphisms continue to hold for quasi-morphisms.

We will denote by $\pp$ the variety consisting of one point. For any variety $X$ let 
$\fL_X=\a_!\bbq\in\cd_mX$ where $\a:X\T T@>>>X$ is the obvious projection. We sometimes write $\fL$ instead 
of $\fL_X$.

Let $v$ be an indeterminate. For any $\ph\in\QQ[v,v\i]$ and any $k\in\ZZ$ we write $(k;\ph)$ for the
coefficient of $v^k$ in $\ph$. Let $\ca=\ZZ[v,v\i]$.

\head Contents\endhead
1. Truncated induction.

2. Truncated restriction. 

3. Truncated convolution from $G_{\e,s}\T G_{\e',s'}$ to $G_{\e\e',s+s'}$.

4. Analysis of the composition $\unz_{\e,s}\unc_{\e,s}$.

5. Adjunction formula (weak form).

6. Equivalence of $\cc^\boc G_{\e,s}$ with the $\e$-centre of $\cc^\boc\cb^2$.

7. Relation with Soergel bimodules.

\head 1. Truncated induction\endhead
\subhead 1.1\endsubhead
For $y\in W$ let $L_y\in\cd_m(\cb^2)$ be the constructible sheaf which is $\bbq$ (with the standard mixed 
structure of pure weight $0$) on $\co_y$ and is $0$ on $\cb^2-\co_y$; let $L_y^\sh\in\cd_m(\cb^2)$ be its
extension to an intersection cohomology complex of $\bco_y$ (equal to $0$ on $\cb^2-\bco_y$). Let 
$\LL_y=L_y^\sh[[|y|+\nu]]\in\cd_m(\cb^2)$.

Let $r\ge1$. For $\ww=(w_1,w_2,\do,w_r)\in W^r$ we set $|\ww|=|w_1|+\do+|w_r|$. Let 
$L^{[1,r]}_\ww\in\cd_m(\cb^{r+1})$ be as in \cite{\TC, 1.1}. For any $J\sub[1,r]$ let
$L_\ww^J\in\cd_m(\cb^{r+1})$, $\dL_\ww^J\in\cd_m(\cb^{r+1})$ be as in \cite{\TC, 1.1}. As in 
\cite{\TC, 1.1(a)}, we have a distinguished triangle
$$(L_\ww^J,L_\ww^{[1,r]},\dL_\ww^J)\tag a$$
in $\cd_m(\cb^{r+1})$. For any $i<i'$ in $[1,r]$ let $p_{i,i'}:\cb^{r+1}@>>>\cb^2$ be the projection to the 
$i,i'$ factors. For ${}^1L,{}^2L,\do,{}^rL$ in $\cd_m(\cb^2)$ we set
$${}^1L\cir{}^2L\cir\do\cir{}^rL=p_{0r!}(p_{01}^*{}^1L\ot p_{12}^*{}^2L\ot\do\ot p_{r-1,r}^*{}^rL)
\in\cd_m(\cb^2).$$ 

\subhead 1.2\endsubhead
Let $\HH$ be the free $\ca$-module with basis $\{T_w;w\in W\}$. It is well known that $\HH$ has
a unique structure of associative $\ca$-algebra with $1=T_1$ (Hecke algebra) such that 
$T_wT_{w'}=T_{ww'}$ if $w,w'\in W$, $|ww'|=|w|+|w'|$ and $T_s^2=1+(v-v\i)T_s$ if $s\in W$, $|s|=1$.
Let $\{c_w;w\in W\}$ be the ``new'' basis of $\HH$ defined as in \cite{\HEC, 5.2} with $L(w)=|w|$.

For $x,y\in W$, the relations $x\prq y$, $x\si y$, $x\si_Ly$ on $W$ are defined as in \cite{\TC, 1.3}. If 
$\boc$ is a two-sided cell of $W$ and $w\in W$, the relations $w\prq\boc$, $\boc\prq w$, $w\prec\boc$, 
$\boc\prec w$ are defined as in \cite{\TC, 1.3}. If $\boc,\boc'$ are two-sided cells of $W$, the relations 
$\boc\prq\boc'$, $\boc\prec\boc'$ are defined as in \cite{\TC, 1.3}. Let $\aa:W@>>>\NN$ be the 
$\aa$-function in \cite{\HEC, 13.6}. If $\boc$ is a two-sided cell of $W$, then for all $w\in\boc$ we have 
$\aa(w)=\aa(\boc)$ where $\aa(\boc)$ is a constant.

Let $\JJ$ be the free $\ZZ$-module with basis $\{t_z;z\in W\}$ with the structure of
associative ring (with $1$) as in \cite{\TC, 1.3}. 
For a two-sided cell $\boc$ of $W$ let $\JJ^{\boc}$ be the subgroup of $\JJ$ generated by 
$\{t_z;z\in\boc\}$; it is a subring of $\JJ$
with unit element $\sum_{d\in\DD_\boc}t_d$ where $\DD_{\boc}$ is the set of distinguished involutions of
$\boc$. We have $\JJ=\op_{\boc}\JJ^\boc$ as rings. 

For $E\in\Irr W$ we define a simple $\QQ\ot\JJ$-module $E_\iy$ and a simple $\QQ(v)\ot_\ca\HH$-module $E(v)$
as in \cite{\TC, 1.3}; there is a unique two-sided cell $\boc_E$ of $W$ such that $\JJ^{\boc_E}E_\iy\ne0$. 

Let $\e\in fA$. Let $E\in\Irr W$. We say that $E\in\Irr_\e W$ if $\tr(\e(w),E)=\tr(w,E)$ for any $w\in W$.
In this case there exists a linear transformation of finite order $\e_E:E@>>>E$
such that $e_Ewe_E\i=\e(w):E@>>>E$ for any $w\in W$; moreover $e_E$ is unique up to multiplication by $-1$.
See \cite{\ORA, 3.2}). For each $E\in\Irr_\e W$ we choose $e_E$ as above.
As a $\QQ$-vector space we have $E_\iy=E$, $E(v)=\QQ(v)\ot_\QQ E$ hence, if
$E\in\Irr_\e W$, $\te:E@>>>E$ can be viewed as a $\QQ$-linear map (of finite order) $\te:E_\iy@>>>E_\iy$ and
as a $\QQ(v)$-linear map (of finite order) $\te:E(v)@>>>E(v)$. From the definitions we see that 
$\te t_w\te\i=t_{\e(w)}:E_\iy@>>>E_\iy$ and $\te T_w\te\i=T_{\e(w)}:E(v)@>>>E(v)$ for any $w\in W$. 

If $E\in\Irr_\e W$ then $\e(\boc_E)=\boc_E$. Let $\Irr_{\e,\boc}W=\{E\in\Irr_\e W;\boc_E=\boc\}$.

\subhead 1.3\endsubhead
For any $\e\in\fA,s\in\ZZ$ let $G_{\e,s}$ be the set of bijections $F:G@>>>G$ such that 

(i) if $s>0$ then $F$ is the Frobenius map for an $F_{p^s}$-rational structure on $G$;

(ii) if $s<0$ then $F\i$ is the Frobenius map for an $F_{p^{-s}}$-rational structure on $G$;

(iii) if $s=0$ then $F$ is an automorphism of $G$;
\nl
moreover in each case (i)--(iii) we require that the following holds: for any $w\in W$ and any 
$(B,B')\in\co_w$ we have $(F(B),F(B'))\in\co_{\e(w)}$.
\nl
(If $\e=1,s=0$ we can identify $G$ and $G_{\e,s}$ by $g\m\Ad(g)$.)
Now $G$ acts on $G_{\e,s}$ by $g:F\m\Ad(g)F\Ad(g\i)$. If $s\ne0$, this action is transitive and the 
stabilizer of a point $F\in G_{\e,s}$ is the finite group $G^F=\{g\in G;F(g)=g\}$. For any $s\in\ZZ$ and any
$\tF\in G_{\e,s}$, the maps $\l:G@>>>G_{\e,s}$, $g\m\Ad(g)\tF$ and $\l':G@>>>G_{\e,s}$, $g\m\tF\Ad(g)$ are 
bijections (by Lang's theorem); we use $\l$ (resp. $\l'$) to view $G_{\e,s}$ with $s\ge0$ (resp. $s\le0$) as
an affine algebraic variety isomorphic to $G$; this algebraic variety 
structure on $G_{\e,s}$ is independent of the choice of $\tF$. We have $\dim G_{\e,s}=\D$.
The $G$-action above on $G_{\e,s}$ is an algebraic group action.
When $X=G_{\e,s}$ then the subscript ${}_m$ in $\cd_m(X),\cm_m(X)$ refers to the 
$\FF_{p^{s_0}}$-structure with Frobenius map $F\m F_0^{s_0}FF_0^{-s_0}$ (with $F_0,s_0$ as in 0.3).

Note that $\sqc_{\e\in\fA,s\in\ZZ}G_{\e,s}$ is a group under composition of maps: if 
$F\in G_{\e,s},F'\in G_{\e',s'}$ then $FF'\in G_{\e\e',s+s'}$. (It is enough to show that for some
$F\in G_{\e,s},F'\in G_{\e',s'}$ we have $FF'\in G_{\e\e',s+s'}$. We take $F=\Ad(\g)F_0^s$, 
$F'=\Ad(\g')F_0^{s'}$ where $\g\in G_{\e,1}$ and $\g'\in G_{\e',1}$ commute with $F_0$; then
$FF'=\Ad(\g\g')F_0^{s+s'}$ and $\g\g'\in G_{\e\e',1}$ commutes with $F_0$ hence
$FF'\in G_{\e\e',s+s'}$.) Note that the composition $G_{\e,s}\T G_{\e',s'}@>>>G_{\e\e',s+s'}$ is not in 
general a morphism of algebraic varieties but only a quasi-morphism (see 0.3), which is good enough for our 
purposes.

{\it Until the end of Section 2 we fix $\e\in\fA$.}
\nl
Let $s\in\ZZ$. We consider the maps $\cb^2@<f<<X_{\e,s}@>\p>>G_{\e,s}$ where 
$$\align&X_{\e,s}=\{(B,B',F)\in\cb\T\cb\T G_{\e,s};F(B)=B'\},\\&f(B,B',F)=(B,B'),\p(B,B',F)=F.\endalign$$ 
Now $L\m\c_{\e,s}(L)=\p_!f^*L$ defines a functor $\cd_m(\cb^2)@>>>\cd_m(G_{\e,s})$.
(When $\e=1,s=0$, $c_{\e,s}$ coincides with the functor $\c$ defined in \cite{\TC, 1.5}). 
For $i\in\ZZ,L\in\cd_m(\cb^2)$ we write $\c_{\e,s}^i(L)$ instead of $(\c_{\e,s}(L))^i$.
For any $z\in W$ we set $R_{\e,s,z}=\c_{\e,s}(L_z^\sh)\in\cd_m(G_{\e,s})$. 
(When $\e=1,s=0$ this is the same as $R_z$ in \cite{\TC, 1.5}.)

Let $b:G_{\e\i,-s}@>\si>>G_{\e,s}$ be the isomorphism $F\m F\i$ and let $b':\cb^2@>\si>>\cb^2$ be the 
isomorphism $(B,B')\m(B',B)$. From the definitions we see that for $L\in\cd_m(\cb^2)$ we have 
$\c_{\e,s}(b'_!L)=b_!\c_{\e\i,-s}(L)$.

Let $CS(G_{\e,s})$ be a set of representatives for the isomorphism classes of simple perverse sheaves 
$A\in\cm(G_{\e,s})$ such that $(A:R_{\e,s,z}^j)\ne0$ for some $z\in W,j\in\ZZ$. (When $\e=1,s=0$ this agrees
with the definition of $CS(G)$ in \cite{\TC, 1.5}.)
Now let $A\in CS(G_{\e,s})$. We associate to $A$ a two-sided cell $\boc_A$ as follows.

Assume first that $s\ne0$. Since $A$ is $G$-equivariant and the conjugation action of $G$ on $G_{\e,s}$ is 
transitive, for any $F\in G_{\e,s}$ we have $A|_{\{F\}}=r_{A,F}[\D]$ where $r_{A,F}$ is an irreducible 
$G^F$-module. From the definitions, for any $z\in W$ and any $F\in G_{\e,s}$ we have 
$$(A:R_{s,z}^j)=(r_{A,F}:IH^{j-\D}\{(B;(B,FB)\in\bco_z\})_{G^F}$$
where the right hand side is the multiplicity of $r_{A,F}$ in the $G^F$-module
$$IH^{j-\D}\{(B;(B,FB)\in\bco_z\};$$
here $IH$ denotes intersection cohomology with coefficients in $\bbq$. In
particular, $r_{A,F}$ is a unipotent representation of $G^F$.
By \cite{\ORA, 3.8}, for any $A\in CS(G_{\e,s})$, any $F\in G_{\e,s}$, any $z\in W$ and any $j\in\ZZ$ we have
$$\align&(r_{A,F}:IH^{j-\D}\{(B;(B,FB)\in\bco_z\})_{G^F}\\&=
(j-\D-|z|;(-1)^{j-\D}\sum_{E\in\Irr_\e W}c_{A,E,\te}\tr(\te c_z,E(v)))\endalign$$
or equivalently
$$(A:R_{\e,s,z}^j)=(j-\D-|z|;(-1)^{j-\D}\sum_{E\in\Irr_\e W}c_{A,E,\te}\tr(\te c_z,E(v)))\tag a$$
where $c_{A,E,\e}$ are uniquely defined rational numbers; now (a) also holds when $s=0$, see
\cite{\TC, 1.5(a)} when $\e=1$ and \cite{\CDGVII, 34.19, 35.22}, \cite{\CDGX, 44.7(e)} for general $\e$.
Moreover, if $s\ne0$ then, by \cite{\ORA, 6.17}, given $A$ as above, there is a unique 
two-sided cell $\boc_A$ of $W$ such that $\e(\boc_A)=\boc_A$ and
$c_{A,E,\e}=0$ whenever $E\in\Irr_\e W$ satisfies $\boc_E\ne\boc_A$. The same holds when $s=0$, see 
\cite{\TC, 1.5} when $\e=1$ and \cite{\CDGIX, \S41} for general $\e$.

When $s\ne0$, $\boc_A$ differs from the two-sided cell associated to $r_{A,F}$ in \cite{\ORA, 4.23} by 
multiplication on the left or right by $w_{max}$. Similarly, when $s=0$, $\boc_A$ differs from the two-sided
cell associated to $A$ in \cite{\CDGIX, \S41} by multiplication on the left or right by $w_{max}$. 

As in \cite{\TC, 1.5(b)}, for $s\in\ZZ$ we have 

(b) {\it $(A:R_{\e,s,z}^j)\ne0$ for some $z\in\boc_A,j\in\ZZ$ and conversely, if $(A:R_{\e,s,z}^j)\ne0$ for 
$z\in W,j\in\ZZ$, then $\boc_A\prq z$.}
\nl
For $s\in\ZZ$, $A\in CS(G_{\e,s})$ let $a_A$ be the value of the $\aa$-function on $\boc_A$. If
$z\in W,E\in\Irr_\e W$ satisfy $\tr(\te c_z,E(v))\ne0$ then $\boc_E\prq z$; if in addition we have 
$z\in\boc_E$, then 
$$\tr(\te c_z,E(v))=\g_{z,E,\te}v^{a_E}+\text{lower powers of }v$$ 
where $\g_{z,E,\e}\in\ZZ$ and $a_E$ is the value of the $\aa$-function on $\boc_E$. Hence from (a) we see 
that

(c) {\it $(A:R_{\e,s,z}^j)=0$ unless $\boc_A\prq z$ and, if $z\in\boc_A$, then}
$$\align&(A:R_{\e,s,z}^j)\\&=(-1)^{j+\D}(j-\D-|z|;(\sum_{E\in\Irr_\e W;\boc_E=\boc_A}c_{A,E,\te}\g_{z,E,\te})
v^{a_A}+\text{ lower powers of }v))\endalign$$
{\it which is $0$ unless $j-\D-|z|\le a_A$.}

{\it In the remainder of this section we fix a two-sided cell $\boc$ of $W$  such that $\e(\boc)=\boc$; we 
set $a=\aa(\boc)$.}
\nl
For $s\in\ZZ$ and $Y=G_{\e,s}$ or $Y=\cb^2$ let $\cm^\spa Y$ be the category of perverse sheaves on $Y$ 
whose composition factors are all of the form $A\in CS(G_{\e,s})$, when $Y=G_{\e,s}$, or of the form $\LL_z$
with $z\in W$, when $Y=\cb^2$. Let $\cm^\prq Y$ (resp. $\cm^\prec Y$) be the category of perverse sheaves on
$Y$ whose composition factors are all of the form $A\in CS(G_{\e,s})$ with $\boc_A\prq\boc$ (resp. 
$\boc_A\prec\boc$), when $Y=G_{\e,s}$, or of the form $\LL_z$ with $z\prq\boc$ (resp. $z\prec\boc$) when 
$Y=\cb^2$. Let $\cd^\spa Y$ (resp. $\cd^\prq Y$ or $\cd^\prec Y$) be the category of all $K\in\cd(Y)$ such 
that $K^i\in\cm^\spa Y$ (resp. $K^i\in\cm^\prec Y$ or $K^i\in\cm^\prec Y$) for all $i\in\ZZ$. Let 
$\cm_m^\spa Y$ (or $\cm^\prq_mY$, or $\cm^\prec_mY$) be the category of all $K\in\cm_mY$ which are also in 
$\cm^\spa Y$ (or $\cm^\prq Y$ or $\cm^\prec Y$). Let $\cd_m^\spa Y$ (or $\cd^\prq_mY$, or $\cd^\prec_mY$) be
the category of all $K\in\cd_mY$ which are also in $\cd^\spa Y$ (or $\cd^\prq Y$ or $\cd^\prec Y$). From (c)
we deduce:

(d) {\it If $z\prq\boc$, then $R_{\e,s,z}^j\in\cm^\prq G_{\e,s}$ for all $j\in\ZZ$. If $z\in\boc$ and 
$j>a+\D+|z|$, then $R_{\e,s,z}^j\in\cm^\prec G_{\e,s}$. If $z\prec\boc$ then 
$R_{\e,s,z}^j\in\cm^\prec G_{\e,s}$ for all $j\in\ZZ$.}

\proclaim{Lemma 1.4}Let $s\in\ZZ$. Let $r\ge1$, $J\sub[1,r]$, $J\ne\emp$ and $\ww=(w_1,w_2,\do,w_r)\in W^r$.
Let $\fE=\D+ra$.

(a) Assume that $w_i\in\boc$ for some $i\in[1,r]$. If $j\in\ZZ$ (resp. $j>\fE$) then 
$\c_{\e,s}^j(p_{0r!}L_\ww^{[1,r]}[|\ww|])$ is in $\cm^\prq G_{\e,s}$ (resp. $\cm^\prec G_{\e,s}$).   

(b) Assume that $w_i\in\boc$ for some $i\in J$. If $j\in\ZZ$ (resp. $j\ge\fE$) then 
$\c_{\e,s}^j(p_{0r!}\dL_\ww^J[|\ww|])$ is in $\cm^\prq G_{\e,s}$ (resp. $\cm^\prec G_{\e,s}$).

(c) Assume that $w_i\in\boc$ for some $i\in J$. If $j\ge\fE$ then the cokernel of the map 
$$\c_{\e,s}^j(p_{0r!}L_\ww^J[|\ww|])@>>>\c_{\e,s}^j(p_{0r!}L_\ww^{[1,r]}[|\ww|])$$
associated to 1.1(a) is in $\cm^\prec G_{\e,s}$.

(d) Assume that $w_i\in\boc$ for some $i\in J$. If $j\in\ZZ$ (resp. $j>\fE$) then 
$\c_{\e,s}^j(p_{0r!}L_\ww^J[|\ww|])$ is in $\cm^\prq G_{\e,s}$ (resp. $\cm^\prec G_{\e,s}$).

(e) Assume that $w_i\prec\boc$ for some $i\in J$. If $j\in\ZZ$ then 
$\c_{\e,s}^j(p_{0r!}L_\ww^{[1,r]}[|\ww|])\in\cm^\prec G_{\e,s}$ and 
$\c_{\e,s}^j(p_{0r!}L_\ww^J[|\ww|])\in\cm^\prec G_{\e,s}$.
\endproclaim  
When $\e=1,s=0$ this is just \cite{\TC, 1.6}; the proof in the general case is entirely similar (it uses 
1.3(b), 1.3(c)).

\subhead 1.5\endsubhead
Let $s\in\ZZ$. Let $CS_{\e,s,\boc}=\{A\in CS(G_{\e,s});\boc_A=\boc\}$. For any $z\in\boc$ we set 
$n_z=a+\D+|z|$. Let $A\in CS_{\e,s,\boc}$ and let $z\in\boc$. We have:
$$(A:R_{s,z}^{n_z})=(-1)^{a+|z|}\sum_{E\in\Irr_{\e,\boc}W}c_{A,E,\te}\tr(\te t_z,E_\iy).\tag a$$
When $\e=1,s=0$ this is just \cite{\TC, 1.7(a)}. In the general case, from 1.3(a) we have
$$(A:R_{s,z}^{n_z})=(-1)^{a+|z|}\sum_{E\in\Irr_\e W}c_{A,E,\te}(a;\tr(\te c_z,E(v)))$$
and it remains to use that $(a;\tr(\te c_z,E(v))$ is equal to $\tr(\te t_z,E_\iy)$ if $E\in\Irr_{\e,\boc}W$
and to $0$, otherwise. We have:

(b) {\it For any $A\in CS_{\e,s,\boc}$ there exists $z\in\boc$ such that $(A:R_{\e,s,z}^{n_z})\ne0$.}
\nl
The proof, based on (a), is the same as that in the case $\e=1,s=0$ given in \cite{\TC, 1.7(b)}.

\mpb

Let $\boc^0=\{z\in\boc;z\si_L\e(z\i)\}$. If $z\in\boc-\boc^0$ and $E\in\Irr_{\e,\boc}W$, then 
$\tr(\te t_z,E_\iy)=0$. (We can write $E_\iy=\op_{d\in\DD_{\boc}}t_dE_\iy$ and $\te t_z:E_\iy@>>>E_\iy$ maps
the summand $t_dE_\iy$ (where $z\si_Ld$) into $t_{\e(d')}E_\iy$ (where $d'\in\DD_\boc$, $d'\si_Lz\i$) and
all other summands to $0$. If $\tr(\te t_z,E_\iy)\ne0$, we must have $t_dE_\iy=t_{\e(d')}E_\iy\ne0$ and 
$d=\e(d')$ and $z\si_L\e(z\i)$.) From this and (a) we deduce

(c) {\it If $z\in\boc-\boc^0$, then $R_{\e,s,z}^{n_z}=0$.}

\subhead 1.6\endsubhead
Let $s\in\ZZ$. For $Y=G_{\e,s}$ or $\cb^2$ let $\cc^\spa Y$ be the subcategory of $\cm^\spa Y$ consisting of 
semisimple objects; let $\cc^\spa_0Y$ be the subcategory of $\cm_mY$ consisting of those $K\in\cm_mY$ such 
that $K$ is pure of weight $0$ and such that as an object of $\cm(Y)$, $K$ belongs to $\cc^\spa Y$. Let 
$\cc^\boc Y$ be the subcategory of $\cm^\spa Y$ consisting of objects which are direct sums of objects in 
$CS_{\e,s,\boc}$ (if $Y=G_{\e,s}$) or of the form $\LL_z$ with $z\in\boc$ (if $Y=\cb^2$). Let $\cc^\boc_0Y$ 
be the subcategory of $\cc^\spa_0Y$ consisting of those $K\in\cc^\spa_0Y$ such that as an object of 
$\cc^\spa Y$, $K$ belongs to $\cc^\boc Y$. For $K\in\cc^\spa_0Y$, let $\un{K}$ be the largest subobject of 
$K$ such that, as an object of $\cc^\spa Y$, we have $\un{K}\in\cc^\boc Y$.   

For $L\in\cc^\spa_0\cb^2$ we define ${}^\e L\in\cc^\spa_0\cb^2$ as follows. We have canonically 
$L=\op_{y\in W}V_y\ot\LL_y$ where $V_y$ are finite dimensional $\bbq$-vector spaces; we set
${}^\e L=\op_{y\in W}V_y\ot\LL_{\e\i(y)}$. We show: 

(a) {\it Let $s\in\NN$. Define $u:G_{\e,s}\T\cb^2@>>>G_{\e,s}\T\cb^2$ by 
$$(F,(B_1,B_2))\m(F,F(B_1),F(B_2))$$ 
and let $L\in\cc^\spa_0\cb^2$. We have canonically $u^*(\bbq\bxt L)=\bbq\bxt{}^\e L$.}
\nl
We can assume that $L=\LL_y$ where $y\in W$; we must show that $u^*(\bbq\bxt\LL_y)=\bbq\bxt\LL_{\e\i(y)}$ or 
that $u^*(\bbq\bxt L^\sh_y)=\bbq\bxt L_{\e\i(y)}^\sh$. Now $\bbq\bxt L^\sh_y$ is the intersection cohomology 
complex of $G_{\e,s}\T\bco_y$ with coefficients in $\bbq$ (extended by $0$ on $G_{\e,s}\T(\cb^2-\bco_y)$).
Hence $u^*(\bbq\bxt L^\sh_y)$ is the intersection cohomology complex of
$u\i(G_{\e,s}\T\bco_y)$ with coefficients in $\bbq$ (extended by $0$ on $G_{\e,s}\T u\i(\cb^2-\bco_y)$)
that is, the intersection cohomology complex of $G_{\e,s}\T\bco_{\e\i(y)}$ with coefficients in $\bbq$ 
(extended by $0$ on $G_{\e,s}\T(\cb^2-\bco_{\e\i(y)})$). This is $\bbq\bxt L^\sh_{\e\i(y)}$, as required.

\mpb

Assume that $s\in\ZZ_{>0}$ and let $F\in G_{\e,s}$. For any $A\in\cc^\spa G_{\e,s}$ we have
$A|_{\{F\}}=r_{A,F}[\D]$ where $r_{A,F}\in\Rep^\spa(G^F)$ (see 0.1). Moreover, from the definitions we see
that

(b) {\it $A\m r_{A,F}$ is an equivalence of categories $\cc^\boc G_{\e,s}@>\si>>\Rep^\boc(G^F)$ (see 0.1).}

\proclaim{Proposition 1.7} Let $s\in\ZZ$.

(a) If $L\in\cd^\prq\cb^2$ then $\c_{\e,s}(L)\in\cd^\prq G_{\e,s}$. If $L\in\cd^\prec\cb^2$, then 
$\c_{\e,s}(L)\in\cd^\prec G_{\e,s}$.

(b) If $L\in\cm^\prq\cb^2$ and $j>a+\nu+\r$ then $\c_{\e,s}^j(L)\in\cm^\prec G_{\e,s}$.
\endproclaim
When $\e=1,s=0$ this is just \cite{\TC, 1.9}; the proof in the general case is entirely similar (it uses 
1.4(a),(e)).

\subhead 1.8\endsubhead
Let $s\in\ZZ$. For $L\in\cc^\boc_0\cb^2$ we set 
$$\unc_{\e,s}(L)=\un{(\c_{\e,s}^{a+\nu+\r}(L)}((a+\nu+\r)/2)
=\un{(\c_{\e,s}(L))^{\{a+\nu+\r\}}}\in\cc^\boc_0G_{\e,s}.$$
The functor $\unc_{\e,s}:\cc^\boc_0\cb^2@>>>\cc^\boc_0G_{\e,s}$ is called {\it truncated induction}. For 
$z\in\boc$ we have 
$$\unc_{\e,s}(\LL_z)=\un{R_{\e,s,z}^{n_z}}(n_z/2).\tag a$$
When $\e=1,s=0$ this is just \cite{\TC, 1.10(a)}; the proof in the general case is entirely similar.

We shall denote by $\t:\JJ^\boc@>>>\ZZ$ the group homomorphism such that $\t(t_z)=1$ if $z\in\DD_\boc$ and 
$\t(t_z)=0$, otherwise. For $z,u\in\boc$ we have:
$$\dim\Hom_{\cc^\boc G_{\e,s}}(\unc_{\e,s}(\LL_z),\unc_{\e,s}(\LL_u))
=\sum_{y\in\boc}\t(t_{y\i}t_zt_{\e(y)}t_{u\i}).\tag b$$
When $\e=1,s=0$ this is just \cite{\TC, 1.10(b)}. We now consider the general case. 

Using (a) and the definitions we see that the left hand side of (b) equals
$$\sum_{A\in CS_{\e,s,\boc}}(A:R_{\e,s,z}^{n_z})(A:R_{\e,s,u}^{n_u}),$$
hence, using 1.5(a) it equals
$$\sum_{E,E'\in\Irr_{\e,\boc}W}
(-1)^{|z|+|u|}\sum_{A\in CS_{\e,s,\boc}}c_{A,E,\te}c_{A,E',\te}\tr(\te t_z,E_\iy)\tr(\te t_u,E'_\iy).$$
Replacing in the last sum $\sum_{A\in CS_{\e,s,\boc}}c_{A,E,\te}c_{A,E',\te}$ by $1$ if $E=E'$ and by $0$ if 
$E\ne E'$ (see \cite{\ORA, 3.9} in the case $s\ne0$ and \cite{\CSIII, 13.12}, \cite{\CDGVII, 35.18(g)} in 
the case $s=0$) we obtain 
$$\sum_{E\in\Irr_{\e,\boc}W}(-1)^{|z|+|u|}\tr(\te t_z,E_\iy)\tr(\te t_u,E_\iy).$$
This is equal to $(-1)^{|z|+|u|}$ times the trace of the operator $\x\m t_z\e(\x)t_{u\i}$ on 
$\QQ\ot\JJ^\boc$ (see \cite{\CDGVII, 34.14(a), 34.17}).
The last trace is equal to the sum over $y\in\boc$ of the coefficient of $t_y$ in $t_zt_{\e(y)}t_{u\i}$; 
this coefficient is equal to $\t(t_{y\i}t_zt_{\e(y)}t_{u\i})$ since for $y,y'\in\boc$, 
$\t(t_{y'}t_y)$ is $1$ if $y'=y\i$ and is $0$ if $y'\ne y\i$ (see \cite{\HEC, 20.1(b)}). Thus we have
$$\dim\Hom_{\cc^\boc G_{\e,s}}(\unc_{\e,s}(\LL_z),\unc_{\e,s}(\LL_u))
=(-1)^{|u|+|z|}\sum_{y\in\boc}\t(t_{y\i}t_zt_{\e(y)}t_{u\i}).$$
Since $\dim\Hom_{\cc^\boc G_{\e,s}}(\unc_{\e,s}(\LL_z),\unc_{\e,s}(\LL_u))\in\NN$ and
$\sum_{y\in\boc}\t(t_{y\i}t_zt_{\e(y)}t_{u\i})\in\NN$, it follows that (b) holds.

\proclaim{Lemma 1.9} Let $s\in\NN$. Let $Y_1,Y_2$ be among $G_{\e,s},\cb^2$ and let $\XX\in\cd_m^\prq Y_1$. 
Let $c,c'$ be integers and let $\Ph:\cd_m^\prq Y_1@>>>\cd_m^\prq Y_2$ be a functor which takes distinguished 
triangles to 
distinguished triangles, commutes with shifts, maps $\cd^\prec_mY_1$ into $\cd^\prec_mY_2$ and maps 
complexes of weight $\le i$ to complexes of weight $\le i$ (for any $i$). Assume that (a),(b) below hold:
$$(\Ph(\XX_0))^h\in\cm_m^\prec Y_2\text{ for any }\XX_0\in\cm_m^\prq Y_1\text{ and any }h>c;\tag a$$
$$\XX\text{ has weight $\le0$ and }\XX^i\in\cm^\prec Y_1\text{ for any }i>c'.\tag b$$
Then 
$$(\Ph(\XX))^j\in\cm^\prec Y_2\text{ for any }j>c+c',\tag c$$
and we have canonically
$$\un{(\Ph(\un{\XX^{\{c'\}}}))^{\{c\}}}=\un{(\Ph(\XX))^{\{c+c'\}}}.\tag d$$
\endproclaim
When $\e=1,s=0$ this is just \cite{\TC, 1.12}; the proof in the general case is entirely similar.

\subhead 1.10\endsubhead
Let $s\in\ZZ$. Let $L\in\cc^\boc_0\cb^2$. We have $\fD(L)\in\cc^\boc_0\cb^2$. Moreover we have canonically:
$$\unc_{\e,s}(\fD(L))=\fD(\unc_{\e,s}(L)).\tag a$$
When $\e=1,s=0$ this is just \cite{\TC, 1.13}; the proof in the general case is entirely similar.

\head 2. Truncated restriction \endhead
\subhead 2.1\endsubhead
Recall that $\e\in\fA$ is fixed.
In this section we fix $s\in\ZZ$. Let $\p,f$ be as in 1.3. Now $K\m\z_{\e,s}(K)=f_!\p^*K$ defines a functor 
$\cd_m(G_{\e,s})@>>>\cd_m(\cb^2)$. (When $\e=1,s=0$, $\z_{\e,s}$ is the same as $\z$ of \cite{\TC, 2.5}.) 
For $i\in\ZZ,K\in\cd_m(G_{\e,s})$ we write $\z_{\e,s}^i(K)$ instead of $(\z_{\e,s}(K))^i$.

Let $b:G_{\e\i,-s}@>\si>>G_{\e,s}$, $b':\cb^2@>\si>>\cb^2$ be as in 1.3. From the definitions we see that 
for $K\in\cd_m(G_{\e\i,-s})$ we have 
$$\z_{\e,s}(b_!K)=b'_!\z_{\e\i,-s}(K).\tag a$$

\proclaim{Proposition 2.2} For any $L\in\cd_m(\cb^2)$ we have
$$\z_{\e,s}(\c_{\e,s}(L))\Bpq
\{\op_{y\in W;|y|=k}L_y\cir L\cir L_{\e(y)\i}\ot\fL[[2k-2\nu]];k\in\NN\},\tag a$$
$$\align&\z_{\e,s}(\c_{\e,s}(L))\Bpq
\\&\{\op_{y\in W;|y|=k}L_y\cir L\cir L_{\e(y)\i}\ot\fL[[2k-2\nu-2\r]]\ot\L^d\cx[[d]](d/2);
k\in\NN,d\in[0,\r]\},\tag b\endalign$$
where $\fL,\cx$ are as in 0.3.
\endproclaim
When $\e=1,s=0$ this is proved in \cite{\TC, 2.6}. 
The proof in the general case will be quite similar to that in the case $\e=1,s=0$. Let
$$Y=\{(B_1,B_2,B_3,B_4,F)\in\cb\T\cb\T\cb\T\cb\T G_s;F(B_1)=B_4,F(B_2)=B_3\}.$$
For $ij=14$ or $23$ we define $h'_{ij}:Y@>>>X_{\e,s}$ by $(B_1,B_2,B_3,B_4,F)\m(B_i,B_j,F)$ and 
$h_{ij}:Y@>>>\cb^2$ by $(B_1,B_2,B_3,B_4,F)\m(B_i,B_j)$. We have $\p^*\p_!=h'_{14!}h'_{23}{}^*$ hence 
$$\z_{\e,s}(\c_{\e,s}(L))=f_!\p^*\p_!f^*(L)=f_!h'_{14!}h'_{23}{}^*f^*(L)=h_{14!}h_{23}^*L.$$
For $k\in\NN$ let $Y^k=\cup_{y\in W;|y|=k}Y_y$ where
$$Y_y=\{(B_1,B_2,B_3,B_4,F)\in Y;(B_1,B_2)\in\co_y,(B_3,B_4)\in\co_{\e(y)\i}\}$$
and let $Y_s^{\ge k}:=\cup_{k';k'\ge k}Y_s^{k'}$, an open subset of $Y_s$; let $h_{ij}^k:Y_s^k@>>>\cb^2$, 
$h_{ij}^{\ge k}:Y_s^{\le k}@>>>\cb^2$ be the restrictions of $h_{ij}$. 
For any $k\in\NN$ we have a distinguished triangle
$$(h_{14!}^{\ge k+1}h_{23}^{\ge k+1*}L),h_{14!}^{\ge k}h_{23}^{\ge k*}L,h_{14!}^kh_{23}^{k*}L).$$
It follows that we have
$$\z_{\e,s}(\c_{\e,s}(L))\Bpq\{h_{14!}^kh_{23}^{k*}L;k\in\NN\}.$$
For $k\in\NN$ let $Z^k=\cup_{y\in W;|y|=k}Z_y$ where
$$Z_y=\{(B_1,B_2,B_3,B_4)\in\cb^4;(B_1,B_2)\in\co_y,(B_3,B_4)\in\co_{\e(y)\i}\};$$
for $i,j\in[1,4]$ we define $\tih_{ij}^k:Z^k@>>>\cb^2$ and $\tih_{ij}^y:Z_y@>>>\cb^2$
by $(B_1,B_2,B_3,B_4)\m(B_i,B_j)$. We have an obvious morphism $u:Y^k@>>>Z^k$.
The fibre of $u$ at $(B_1,B_2,B_3,B_4)\in Z^k$ can be identified with the set of all
$F\in G_{\e,s}$ such that $F(B_1)=B_4,F(B_2)=B_3$.
Since $(B_1,B_2)\in\co_y,(B_3,B_4)\in\co_{\e(y)\i}$ for some $y\in W$, we can find $\tF\in G_{\e,s}$ such 
that $\tF(B_1)=B_4,\tF(B_2)=B_3$; hence the fibre above can be identified with 
$$\align&\{g\in G;\Ad(g)\tF(B_1)=B_4,\Ad(g)\tF(B_2)=B_3\}\\&
=\{g\in G;\Ad(g)(B_4)=B_4,\Ad(g)(B_3)=B_3\}=B_3\cap B_4\endalign$$ 
which is quasi-isomorphic to $\kk^{\nu-k}$ times the $\r$-dimensional torus $T$. We have a commutative 
diagram
$$\CD
\cb^2@<h_{23}^k<<Y_s^k@>h_{14}^k>>\cb^2\\
@V1VV     @VuVV   @V1VV \\
\cb^2@<\tih_{23}^k<<Z^k@>\tih_{14}^k>>\cb^2\endCD$$
We have
$$h_{14!}^kh_{23}^{k*}L=\tih_{14!}^ku_!u^*\tih_{23}^{k*}L=\tih_{14!}^k(\tih_{23}^{k*}L\ot u_!\bbq)=
(\tih_{14!}^k\tih_{23}^{k*}L)\ot\fL[[-2\nu+2k]].$$  
We deduce that
$$\z_{\e,s}(\c_{\e,s}(L))\Bpq\{(\tih_{14!}^k\tih_{23}^{k*}L)\ot\fL[[-2\nu+2k]];k\in\NN\}.$$
Since $Z^k$ is the union of open and closed subvarieties $Z_y,|y|=k$, we have
$$\tih_{14!}^k\tih_{23}^{k*}L=\op_{y\in W;|y|=k}\tih_{14!}^y\tih_{23}^{y*}L.$$
From the definitions we have
$$\tih_{14!}^y\tih_{23}^{y*}L=L_y\cir L\cir L_{\e(y)\i}.$$
This completes the proof of (a). Now (b) follows from (a) just as in the case where $\e=1,s=0$.

{\it In the remainder of this section we fix a two-sided cell $\boc$ of $W$  such that $\e(\boc)=\boc$; we 
set $a=\aa(\boc)$.}

\proclaim{Proposition 2.3}Let $w\in W$ and let $j\in\ZZ$. We set 
$S=\z_{\e,s}(R_{\e,s,w})[[2\r+2\nu+|w|]]\in\cd_m(\cb^2)$.

(a) If $w\prq\boc$ then $S^j\in\cm^\prq\cb^2$.

(b) If $w\in\boc$ and $j>\nu+2a$ then $S^j\in\cm^\prec\cb^2$.

(c) If $w\prec\boc$ then $S^j\in\cm^\prec\cb^2$.

(d) $S^j$ is mixed of weight $\le j$.

(e) If $j\ne\nu+2a$ and $w\in\boc$ then $gr_{\nu+2a}S^j\in\cm^\prec\cb^2$.

(f) If $k>\nu+2a$ and $w\in\boc$ then $gr_kS^j\in\cm^\prec\cb^2$.
\endproclaim
When $\e=1,s=0$ this is just \cite{\TC, 2.7}. The proof in the general case is entirely similar; it uses 2.2.

\proclaim{Proposition 2.4}(a) If $K\in\cd^\prq G_{\e,s}$ then $\z_{\e,s}(K)\in\cd^\prq\cb^2$. If 
$K\in\cd^\prec G_{\e,s}$, then $\z_{\e,s}(K)\in\cd^\prec\cb^2$.

(b) If $K\in\cm^\prq G_{\e,s}$ and $j>\r+\nu+a$ then $\z_{\e,s}^j(K)\in\cm^\prec\cb^2$.
\endproclaim
When $\e=1,s=0$ this is just \cite{\TC, 2.8}. The proof in the general case is entirely similar; it uses 
1.5(b) and 2.3.

\subhead 2.5\endsubhead
For $K\in\cc^\boc_0G_{\e,s}$ we set 
$$\unz_{\e,s}(K)=\un{(\z_{\e,s}(K))^{\{\r+\nu+a\}}}\in\cc^\boc_0\cb^2.$$
We say that $\unz_{\e,s}(K)$ is the {\it truncated restriction} of $K$.

\proclaim{Proposition 2.6} Let $K\in\cd_m(G_{\e,s})$ and let $L\in\cc^\spa_0\cb^2$. Then there is a 
canonical isomorphism ${}^\e L\cir\z_{\e,s}(K)@>\si>>\z_{\e,s}(K)\cir L$.
\endproclaim                                        
When $\e=1,s=0$ this follows from \cite{\TC, 2.10(a)}. We now consider the general case. Let 
$u:G_{\e,s}\T\cb^2@>>>G_{\e,s}\T\cb^2$ be as in 1.6(a).
We have $\z_{\e,s}(K)\cir L=c_!d^*(K\bxt L)$ where $Z=\{(F,(B,B'',B'))\in G_{\e,s}\T\cb^3;F(B)=B''\},$
$d:Z@>>>G_{\e,s}\T\cb^2$ is $(F,(B,B'',B'))\m(F,(B'',B'))$, $c:Z@>>>\cb^2$ is $(F,(B,B'',B'))\m(B,B')$. We 
have ${}^\e L\cir\z_{\e,s}(K)=c'_!d'{}^*(K\bxt {}^\e L)$ where 
$Z'=\{(F,(B,B'',B'))\in G_{\e,s}\T\cb^3;F(B'')=B'\}$, $d':Z'@>>>G_{\e,s}\T\cb^2$ is 
$(F,(B,B'',B'))\m(F,(B,B''))$, $c':Z'@>>>\cb^2$ is $(F,(B,B'',B'))\m(B,B')$.
Using 1.6(a) we have $K\bxt{}^\e L=u^*(K\bxt L)$ hence it is enough to show that
$c_!d^*(K\bxt L)=c'_!d'{}^*u^*(K\bxt L)$.
We have $c_!d^*(K\bxt L)=c_{1!}d_1^*(K\bxt L)=c'_!d'{}^*u^*(K\bxt L)$  where
$d_1:G_{\e,s}\T\cb^2@>>>G_{\e,s}\T\cb^2$ is $(F,(B,B'))\m(F,(F(B),B'))$, $c_1:G_{\e,s}\T\cb^2@>>>\cb^2$ is
$(F,(B,B'))\m(B,B')$. The proposition follows.

\proclaim{Proposition 2.7} (a) If $L\in\cm^\prq\cb^2$ and $j>2a+2\nu+2\r$ then 
$(\z_{\e,s}(\c_{\e,s}(L)))^j\in\cm^\prec\cb^2$.

(b) If $L\in\cc^\boc_0\cb^2$, we have canonically 
$$\unz_{\e,s}(\unc_{\e,s}(L))=\un{(\z_{\e,s}(\c_{\e,s}(L)))^{\{2a+2\nu+2\r\}}}\in\cc^\boc_0\cb^2.$$
\endproclaim
We apply 1.9 with $\Ph=\z_{\e,s}:\cd_m(G_{\e,s})@>>>\cd_m(\cb^2)$ and with $\XX=\c_{\e,s}(L)$,
$(c,c')=(a+\nu+\r,a+\nu+\r)$, see 2.4, 1.7. The result follows.

\subhead 2.8\endsubhead
For $L,L'\in\cc^\boc_0\cb^2$, we set (as in \cite{\TC, 3.2})
$$L\unb L'=\un{(L\cir L')^{\{a-\nu\}}}\in\cc^\boc_0\cb^2.\tag a$$
This defines an associative tensor product structure on $\cc^\boc_0\cb^2$.

\proclaim{Proposition 2.9}Let $K\in\cc^\boc_0G_{\e,s},L\in\cc^\boc_0\cb^2$. There is a canonical isomorphism
$${}^\e L\unb\unz_{\e,s}(K)@>\si>>\unz_{\e,s}(K)\unb L.\tag a$$
\endproclaim
Applying 1.9 with $\Ph:\cd_m^\prq\cb^2@>>>\cd_m^\prq\cb^2$, $L'\m L'\cir L$, $\XX=\z_{\e,s}(K)$,
$(c,c')=(a-\nu,a+\r+\nu)$ (see \cite{\TC, 3.1} and 2.4), we deduce that we have canonically 
$$\un{((\z_{\e,s}(K))^{\{a+\r+\nu\}}\cir L)^{\{a-\nu\}}}=\un{(\z_{\e,s}(K)\cir L)^{\{2a+\r\}}}.\tag b$$
Using 1.9 with $\Ph:\cd_m^\prq\cb^2@>>>\cd_m^\prq\cb^2$, $L'\m{}^\e L\cir L'$, $\XX=\z_{\e,s}(K)$,
$(c,c')=(a-\nu,a+\r+\nu)$ (see \cite{\TC, 3.1} and 2.8), we deduce that we have canonically 
$$\un{({}^\e L\cir(\z_{\e,s}(K))^{\{a+\r+\nu\}})^{\{a-\nu\}}}=\un{({}^\e L\cir\z_{\e,s}(K))^{\{2a+\r\}}}.
\tag c$$
We now combine (b),(c) with 2.6; we obtain the isomorphism (a).

\subhead 2.10\endsubhead
Define $c:G_{\e,s}\T\cb^2@>>>\cb^2$ by $(F,B,B')\m(F(B),F(B'))$. We show that for $K\in\cc^\spa G_{\e,s}$
we have canonically 
$$c^*\z_{\e,s}K=\bbq\bxt\z_{\e,s}K.\tag a$$
We have a commutative diagram with cartesian left squares
$$\CD
G_{\e,s}\T\cb^2@<f''<< X''_{\e,s}@>\p''>>G_{\e,s}\T G_{\e,s}@>e>>G_{\e,s}\\
@VdVV                  @Vd'VV           @Vd''VV\\     @.
G_{\e,s}\T\cb^2@<f'<< X'_{\e,s}@>\p'>>G_{\e,s}@.       {}
@VcVV                   @Vc'VV        @Vc''VV     @.\\
\cb^2@<f<<X_{\e,s}@>\p>>G_{\e,s} @. {}\\
\endCD$$
where $f,g$ are as in 1.3, 

$X'_{\e,s}=\{(\tF,B,B',F)\in G_{\e,s}\T\cb\T\cb\T G_{\e,s};F\tF(B)=\tF(B')\}$,

$X''_{\e,s}=\{(\tF,B,B',F)\in G_{\e,s}\T\cb\T\cb\T G_{\e,s};F(B)=B'\}$,

$f'(\tF,B,B',F)=(\tF,B,B')$, $f''(\tF,B,B',F)=(\tF,B,B')$, 

$\p'(\tF,B,B',F)=F$, $\p''(\tF,B,B',F)=(\tF,F)$, 

$c'(\tF,B,B',F)=(F,\tF(B),\tF(B'))$, $c''(F)=F$, $d(\tF,B,B')=(\tF,B,B')$,

$d'(\tF,B,B',F)=(\tF,\tF\i F\tF,B,B')$, $d''(\tF,F)=\tF\i F\tF$, $e(\tF,F)=F$ .
\nl
It is enough to show that $d^*f'_!c'{}^*\p^*K=f''_!\p''{}^*e^*K$. or that
$f''_!d'{}^*\p'{}^*K=f''_!\p''{}^*e^*K$. It is enough to show that $d'{}^*\p'{}^*K=\p''{}^*e^*K$, or that 
$\p''{}^*d''{}^*K=\p''{}^*e^*K$. Hence it is enough to show that $d''{}^*K=e^*K$. We identify 
$G\T G_{\e,s}\lra G_{\e,s}\T G_{\e,s}$ by $(g,F)\lra(F\Ad(g),F)$. Then
$d'',e:G_{\e,s}\T G_{\e,s}@>>>G_{\e,s}$ become the maps $d_1,e_1:G\T G_{\e,s}@>>>G_{\e,s}$ given by
$(g,F)\m\Ad(g(\i F\Ad(g)$, $(g,F)=F$ respectively and we have $d_1^*K=e_1^*K$ by the $G$-equivariance of $K$.
Hence $d''{}^*K=e^*K$ as required.

Using (a) and the definitions we see that for any $K\in\cc^\boc_0G_{\e,s}$ we have canonically 
$$c^*\unz_{\e,s}K=\bbq\bxt\unz_{\e,s}K.\tag b$$
From the definitions (see 1.6) for any $L\in\cc^\spa_0\cb^2$ we have $c^*L=\bbq\bxt{}^\e L$. Comparing with
(b) we deduce that we have canonically
$${}\e(\unz_{\e,s}K)=\unz_{\e,s}K\tag c$$
for any $K\in\cc^\boc_0G_{\e,s}$.

\head 3. Truncated convolution from $G_{\e,s}\T G_{\e',s'}$ to $G_{\e\e',s+s'}$\endhead
\subhead 3.1\endsubhead
Let $\e,\e'\in\fA$, $s,s'\in\ZZ$. We define $\mu:G_{\e,s}\T G_{\e',s'}@>>>G_{\e\e',s+s'}$ by $(F,F')=FF'$ 
(composition of maps $G@>>>G$); this is a quasi-morphism, see 1.3. For 
$K\in\cd_m(G_{\e,s}),K'\in\cd_m(G_{\e',s'})$ we define the convolution $K*K'\in\cd_m(G_{\e\e',s+s'})$ by 
$K*K'=\mu_!(K\bxt K')$. 
If $\e''\in\fA$, $s''\in\ZZ$ then for $K,K'$ as above and $K''\in\cd_m(G_{\e'',s''})$, we have canonically 
$(K*K')*K''=K*(K'*K'')\in\cd_m(G_{\e\e'\e'',s+s'+s''})$ (and we denote this by $K*K'*K''$).

\proclaim{Lemma 3.2} Let $\e,\e'\in\fA$, $s,s'\in\ZZ$. Let $K\in\cd_m(G_{\e,s})$, $L\in\cd_m(\cb^2)$. 
We have canonically $K*\c_{\e',s'}(L)=\c_{\e\e',s+s'}(L\cir\z_{\e,s}(K))$.
\endproclaim
Let 
$$Z=\{(F_1,F_2,B,B')\in G_{\e,s}\T G_{\e',s'}\T\cb\T\cb;F_2(B)=B'\}.$$ 
Define $c:Z@>>>G_{\e,s}\T\cb^2$ by
$(F_1,F_2,B,B')\m(F_1,(B,B'))$ and $d:Z@>>>G_{\e\e',s+s'}$ by $(F_1,F_2,B,B')\m F_1F_2$. From the 
definitions we see that both 
$$K*\c_{\e',s'}(L),\c_{\e\e',s+s'}(L\cir\z_{\e,s}(K))$$ 
can be identified with $d_!c^*(K\bxt L)$. The lemma follows. (In the case where $\e=\e'=1$ and $s=s'=0$ this
reduces to \cite{\TC, 4.2}.)

\proclaim{Proposition 3.3} Let $\e,\e'\in\fA$, $s,s'\in\ZZ$. For any $L,L'\in\cd_m(\cb^2)$ we have
$$\align&\c_{\e,s}(L)*\c_{\e',s'}(L')[[2\r+2\nu]]\\&\Bpq  
\{\c_{\e\e',s+s'}(L'\cir L_y\cir L\cir L_{\e(y)\i})[[2|y|]]\ot\L^d\cx[[d]](d/2);d\in[0,\r],y\in W\}.
\endalign$$
\endproclaim
From 2.2(b) we deduce
$$\align&L'\cir\z_{\e,s}(\c_{\e,s}(L))[[2\nu+2\r]]\\&
\Bpq\{L'\cir L_y\cir L\cir L_{\e(y)\i}[[2|y|]]\ot\L^d\cx[[d]](d/2);y\in W,d\in[0,\r]\}\endalign$$
and
$$\align&\c_{\e\e',s+s'}(L'\cir\z_{\e,s}(\c_{\e,s}(L)))[[2\nu+2\r]]\\&
\Bpq\{\c_{\e\e',s+s'}(L'\cir L_y\cir L\cir L_{\e(y)\i})[[2|y|]]\ot\L^d\cx[[d]](d/2);y\in W,d\in[0,\r]\}.
\endalign$$
It remains to show that $\c_{\e\e',s+s'}(L'\cir\z_{\e,s}(\c_{\e,s}(L)))=\c_{\e,s}(L)*\c_{\e',s'}(L')$. This 
follows from 3.2 with $K,L$ replaced by $\c_{\e,s}(L),L'$.

{\it In the remainder of this section we fix a two-sided cell $\boc$ of $W$; we set $a=\aa(\boc)$.}

\proclaim{Proposition 3.4} Let $\e,\e'\in\fA$, $s,s'\in\ZZ$. Assume that $\e(\boc)=\boc$, $\e'(\boc)=\boc$.
Let $w,w'\in W$ and let $j\in\ZZ$. We set 
$C=R_{\e,s,w}*R_{\e',s',w'}[[2\r+2\nu+|w|+|w'|]]\in\cd_m(G_{\e\e',s+s'})$.

(a) If $w\prq\boc$ or $w'\prq\boc$ then $C^j\in\cm^\prq G_{\e\e',s+s'}$.

(b) If $j>\D+4a$ and either $w\in\boc$ or $w'\in\boc$ then $C^j\in\cm^\prec G_{\e\e',s+s'}$.

(c) If $w\prec\boc$ or $w'\prec\boc$ then $C^j\in\cm^\prec G_{\e\e',s+s'}$.

(d) $C^j$ is mixed of weight $\le j$.

(e) If $j\ne\D+4a$ and either $w\in\boc$ or $w'\in\boc$ then $gr_{\D+4a}C^j\in\cm^\prec G_{\e\e',s+s'}$.

(f) If $k>\D+4a$ and $w\in\boc$ or $w'\in\boc$ then $gr_k C^j\in\cm^\prec G_{\e\e',s+s'}$.
\endproclaim
When $\e=\e'=1$, $s=s'=0$, this is just \cite{\TC, 4.4}. The proof in the general case is entirely similar; 
it uses 3.3 and 1.4(d),(e).

\proclaim{Proposition 3.5} Let $\e,\e'\in\fA$, $s,s'\in\ZZ$. Assume that $\e(\boc)=\boc$, $\e'(\boc)=\boc$.
Let $K\in\cd^\spa_m(G_{\e,s})$, $K'\in\cd^\spa_m(G_{\e',s'})$.

(a) If $K\in\cd^\prq G_{\e,s}$ or $K'\in\cd^\prq G_{\e',s'}$ then $K*K'\in\cd^\prq G_{\e\e',s+s'}$; if 
$K\in\cd^\prec G_{\e,s}$ or $K'\in\cd^\prec G_{\e',s'}$ then $K*K'\in\cd^\prec G_{\e\e',s+s'}$.

(b) If $K\in\cm^\prq G_{\e,s}$, $K'\in\cm^\prq G_{\e',s'}$ and $j>\r+2a$ then 
$(K*K')^j\in\cm^\prec G_{\e\e',s+s'}$.
\endproclaim 
When $\e=\e'=1$, $s=s'=0$, this is just \cite{\TC, 4.5}. The proof in the general case is entirely similar; 
it uses 3.4.

\subhead 3.6\endsubhead
Let $\e,\e'\in\fA$, $s,s'\in\ZZ$. Assume that $\e(\boc)=\boc$, $\e'(\boc)=\boc$. For 
$K\in\cc^\boc_0G_{\e,s}$, $K'\in\cc^\boc_0G_{\e',s'}$ we set
$$K\unst K'=\un{(K*K')^{\{2a+\r\}}}\in\cc^\boc_0G_{\e\e',s+s'}.$$
We say that $K\unst K'$ is the {\it truncated convolution} of $K,K'$.

\proclaim{Proposition 3.7} Let $\e,\e',\e''\in\fA$, $s,s',s''\in\ZZ$. 
Assume that $\e(\boc)=\boc$, $\e'(\boc)=\boc$, $\e''(\boc)=\boc$. Let $K,K',K''$ be in 
$\cc^\boc_0G_{\e,s},\cc^\boc_0G_{\e',s'},\cc^\boc_0G_{\e'',s''}$ respectively. There is a canonical 
isomorphism
$$(K\unst K')\unst K''@>\si>>K\unst(K'\unst K'').\tag a$$
\endproclaim
When $\e=\e'=\e''=1$, $s=s'=s''=0$, this is just \cite{\TC, 4.7}. The proof in the general case is entirely 
similar; it uses 1.9, 3.5.

\proclaim{Proposition 3.8}Let $\e,\e'\in\fA$, $s,s'\in\ZZ$. Assume that $\e(\boc)=\boc$, $\e'(\boc)=\boc$. 
Let $K\in\cc^\boc_0G_{\e,s}$, $K'\in\cc^\boc_0G_{\e',s'}$. There is a canonical isomorphism (in 
$\cc^\boc_0\cb^2$):
$$\unz_{\e',s'}(K')\unb\unz_{\e,s}(K)@>\si>>\unz_{\e\e',s+s'}(K\unst K').$$
\endproclaim
When $\e=\e'=1$, $s=s'=0$ this is just \cite{\TC, 5.2}. The proof in the general case is entirely similar.

\head 4. Analysis of the composition $\unz_{\e,s}\unc_{\e,s}$\endhead
\subhead 4.1\endsubhead
In the remainder of this paper we fix a two-sided cell $\boc$ of $W$; we set $a=\aa(\boc)$. We also fix
$\e\in\fA$ such that $\e(\boc)=\boc$. In this section we fix $s\in\ZZ$. 
Let $e,f,e'$ be integers such that $e\le f\le e'-3$ and let $\ee=e'-e+1$; we have $\ee\ge4$. We set 
$$\cy=\{((B_e,B_{e+1},\do,B_{e'}),F)\in\cb^\ee\T G_s;F(B_f)=B_{f+3},F(B_{f+1})=B_{f+2}\}.$$
Define $\vt:\cy@>>>\cb^\ee$ by $((B_e,B_{e+1},\do,B_{e'}),F)\m(B_e,B_{e+1},\do,B_{e'})$. For $i,j$ in 
$\{e,e+1,\do,e'\}$ let $p_{ij}:\cb^e@>>>\cb^2$ be the projection to the $i,j$ coordinate; define 
$h_{ij}:\cy@>>>\cb^2$ by $h_{ij}=p_{ij}\vt$. Now $G^{\ee-2}$ acts on $\cy$ by 
$$\align&(g_e,\do,g_f,g_{f+3},\do,g_{e'}):((B_e,B_{e+1},\do,B_{e'}),F)\m\\&
(\Ad(g_e)(B_e),\Ad(g_{e+1}((B_{e+1}),\do,\Ad(g_{f-1})(B_{f-1}),\Ad(g_f)(B_f),\Ad(g_f)(B_{f+1}),\\&
\Ad(g_{f+3})(B_{f+2}),\Ad(g_{f+3})(B_{f+3}),\Ad(g_{f+4})(B_{f+4}),\do,\Ad(g_{e'})(B_{e'})),\\&
\Ad(g_{f+3})F\Ad(g_f\i));\endalign$$
this induces a $G^{\ee-2}$-action on $\cb^\ee$ so that $\vt$ is $G^{\ee-2}$-equivariant.

Let $E=\{e,e+1,\do,e'-1\}-\{f,f+2\}$. Assume that $x_n\in\boc$ are given for $n\in E$. 
Let $P=\ot_{n\in E}p_{n,n+1}^*\LL_{x_n}\in\cd_m\cb^\ee$,
$\tP=\ot_{n\in E}h_{n,n+1}^*\LL_{x_n}=\vt^*P\in\cd_m\cy$.
In 4.1-4.7 we will study 
$$h_{ee'!}\tP\in\cd_m\cb^2.$$
Setting $\Xi=\vt_!\bbq\in\cd_m\cb^\ee$, we have 
$$h_{ee'!}\tP=p_{ee!}(\Xi\ot P).$$
Clearly, $\Xi^j$ is $G^{\ee-2}$-equivariant for any $j$. For any $y,y'$ in $W$ we set
$$Z_{y,y'}:=\{(B_e,B_{e+1},\do,B_{e'})\in\cb^\ee;(B_f,B_{f+1})\in\co_y,(B_{f+2},B_{f+3})\in\co_{y'}\}.$$
These are the orbits of the $G^{\ee-2}$-action on $\cb^\ee$. Note that the fibre of $\vt$ over a point of 
$Z_{y,y'}$ is isomorphic to $T\T\kk^{\nu-|y|}$ if $y'=\e(y)\i$ and is empty if $y'\ne\e(y)\i$. Thus 

(a) $\Xi|_{Z_{y,y'}}$ is $0$ if $y'\ne\e(y)\i$
\nl
and for any $y\in W$ we have  
$$\ch^h\Xi|_{Z_{y,\e(y)\i}}=0\text{ if }h>2\nu-2|y|+2\r,\qua 
\ch^{2\nu-2|y|+2\r}\Xi|_{Z_{y,\e(y)\i}}=\bbq(-\nu+|y|-\r).\tag b$$
The closure of $Z_{y,y'}$ in $\cb^\ee$ is denoted by $\bZ_{y,y'}$. We set $k_\ee=\ee\nu+2\r$. We have the 
following result.

\proclaim{Lemma 4.2} (a) We have $\Xi^j=0$ for any $j>k_\ee$. Hence, setting $\Xi'=\t_{\le k_\ee-1}\Xi$, we 
have a canonical distinguished triangle $(\Xi',\Xi,\Xi^{k_\ee}[-k_\ee])$. 

(b) If $\x\in Z_{y,y'}$ and $i=2\nu-|y|-|y'|+2\r$, the induced homomorphism 
$\ch^i_\x\Xi@>>>\ch^{i-k_\ee}_\x(\Xi^{k_\ee})$ is an isomorphism.
\endproclaim
When $\e=1,s=0$ this is just \cite{\TC, 6.2}. The proof in the general case is entirely similar; it uses 
4.1(a),(b).

\subhead 4.3\endsubhead
For any $y,y'$ in $W$ let $\fT_{y,y'}$ be the intersection cohomology complex of $\bZ_{y,y'}$ extended by 
$0$ on $\cb^\ee-\bZ_{y,y'}$, to which $[[(\ee-2)\nu+|y|+|y'|]]$ is applied. Note that
$$\fT_{y,y'}=p_{f,f+1}^*\LL_y\ot p_{f+2,f+3}^*\LL_{y'}[[(\ee-4)\nu]].\tag a$$
We have the following result.

\proclaim{Lemma 4.4} We have canonically $gr_0(\Xi^{k_\ee}(k_\ee/2))=\op_{y\in W}\fT_{y,\e(y)\i}$.
\endproclaim
When $\e=1,s=0$ this is just \cite{\TC, 6.4}. The proof in the general case is entirely similar; it uses 
4.2(b) and 4.1.

\subhead 4.5\endsubhead
Let $y,\ty\in W$. Using the definitions and 1.2(a) we have 
$$\align&p_{ee'!}(\fT_{y,\ty}\ot P[[(6-2\ee)\nu]])\\&
=L_{x_1}^\sh\cir\do\cir L_{x_{f-1}}^\sh\cir L_y^\sh\cir L_{x_{f+1}}^\sh\cir L_{\ty}^\sh\cir
L_{x_{f+3}}^\sh\cir\do\cir L_{x_{e'}}^\sh[[\nu+|y|+|\ty|+\sum_{n\in E}|x_n|]].\tag a\endalign$$

\proclaim{Lemma 4.6}The map $\Xi@>>>\Xi^{k_\ee}[-k_\ee]$ (coming from $(\Xi',\Xi,\Xi^{k_\ee}[-k_\ee])$ in 
4.2(a)) induces a morphism 
$$(p_{ee'!}(\Xi\ot P))^{(\ee-2)a+(6-\ee)\nu+2\r}@>>>(p_{ee'!}(\Xi^{k_\ee}\ot P))^{(\ee-2)a+(6-\ee)\nu+2\r
-k_\ee}$$
whose kernel and cokernel are in $\cm^\prec_m\cb^2$.
\endproclaim
When $\e=1,s=0$ this is just \cite{\TC, 6.6}. The proof in the general case is entirely similar; it uses 
4.5(a),(b) and \cite{\TC, 2.2(a)}.

\proclaim{Lemma 4.7} We have canonically
$$\un{(h_{ee'!}\tP)^{\{(\ee-2)a+(6-\ee)\nu+2\r\}}}=\op_{y\in\boc}Q_y$$
where 
$$\align&Q_y=\un{(p_{ee'!}(\fT_{y,\e(y)\i}\ot P))^{\{(\ee-2)a+(6-2\ee)\nu\}}}\\&
=\LL_{x_1}\unb\do\unb\LL_{x_{f-1}}\unb\LL_y\unb\LL_{x_{f+1}}\unb\LL_{\e(y)\i}\unb\LL_{x_{f+3}}\unb\do
\unb\LL_{x_{e'}}.\endalign$$
\endproclaim
When $\e=1,s=0$ this is just \cite{\TC, 6.7}. The proof in the general case is entirely similar; it uses 4.6,
4.5(a),(b) and \cite{\TC, 2.2(a), 2.3, 3.2}.

\proclaim{Theorem 4.8} Let $x\in\boc$. We have canonically
$$\unz_{\e,s}(\unc_{\e,s}(\LL_x))=\op_{y\in\boc}\LL_y\unb\LL_x\unb\LL_{\e(y)\i}.\tag a$$
\endproclaim
When $\e=1,s=0$ this is just \cite{\TC, 6.8}. The proof in the general case is entirely similar; it uses 4.7,
the proof of 2.2 and 2.7(b).

\subhead 4.9\endsubhead
Using \cite{\TC, 2.4} we see that 4.8(a) implies
$$\unz_{\e,s}\unc_{\e,s}\LL_x\cong\op_{z\in\boc}(\LL_z)^{\op\ps_x(z)}\tag a$$
in $\cc^\boc\cb^2$ where $\ps_x(z)\in\NN$ are given by the following equation in $\JJ^\boc$:
$$\sum_{y\in\boc}t_yt_xt_{\e(y)\i}=\sum_{z\in\boc}\ps_x(z)t_z.$$

\head 5. Adjunction formula (weak form)\endhead
\proclaim{Proposition 5.1} Let $\e'\in\fA$, $s,s'\in\ZZ$. We assume that $\e'(\boc)=\boc$.
Let $K\in\cc^\boc_0(G_{\e,s})$, $L\in\cc^\boc_0(\cb^2)$. We have canonically 
$$K\unst\unc_{\e',s'}(L)=\unc_{\e\e',s+s'}(L\unb\unz_{\e,s}(K)).\tag a$$
\endproclaim
When $\e=\e'=1,s=s'=0$ this is just \cite{\TC, 8.1}. The proof in the general case is entirely similar; it 
uses 3.2.

\subhead 5.2\endsubhead
Let $s\in\ZZ$. When $\e=1,s=0$, the arguments in this subsection reduce to arguments in \cite{\TC, 8.8}.
Let $u':G_{\e\i,-s}@>>>\pp$ be the obvious map. From \cite{\CSII, 7.4} we see that for 
$K,K'\in\cm_m^\prq G_{\e\i,-s}$ we have canonically
$$(u'_!(K\ot K'))^0=\Hom_{\cm G_{\e\i,-s}}(\fD(K),K'),\qua (u'_!(K\ot K'))^j=0\text { if }j>0.$$
We deduce that if $K,K'$ are also pure of weight $0$ then $(u'_!(K\ot K'))^0$ is pure of weight zero that 
is, $(u'_!(K\ot K'))^0=gr_0(u'_!(K\ot K'))^0$. Let $\io:\pp@>>>G=G_0$ be the map with image $1$. From the 
definitions we see that we have $u'_!(K\ot K')=\io^*(b_!(K)*K')$ where $b:G_{\e\i,-s}@>>>G_{\e,s}$ is given 
by $F\m F\i$. Hence for $K,K'\in\cc^\boc_0G_{\e\i,-s}$ we have  
$$\Hom_{\cc^\boc G_{\e\i,-s}}(\fD(K),K')=(\io^*(b_!(K)*K'))^0=(\io^*(b_!(K)*K'))^{\{0\}}.\tag a$$
Applying \cite{\TC, 8.2} with $\Ph:\cd_m^\prq G_0@>>>\cd_m\pp$, $K_1\m\io^*K_1$, $c=-2a-\r$ (see 
\cite{\TC, 8.3(a)}), $K$ replaced by $b_!(K)*K'\in\cd_m(G_{1,0}$ and $c'=2a+\r$ we see that we have 
canonically
$$(\io^*(b_!(K)\unst K'))^{\{-2a-\r\}}\sub(\io^*(b_!(K)*K'))^{\{0\}}.$$
In particular, if $L,L'\in\cc^\boc_0\cb^2$ then we have canonically  
$$(\io^*(\unc_{\e,s}(L')\unst\unc_{\e\i,-s}(L)))^{\{-2a-\r\}}\sub
(\io^*(\unc_{\e,s}(L')*\unc_{\e\i,-s}(L)))^{\{0\}}.$$
Using the equality 
$$(\io^*(\unc_{\e,s}(L')\unst\unc_{\e\i,-s}(L)))^{\{-2a-\r\}}
=(\io^*(\unc_{1,0}(L\unb\unz_{\e,s}(\unc_{\e,s}(L')))))^{\{-2a-\r\}}$$
which comes from 5.1, we deduce that we have canonically
$$(\io^*(\unc_{1,0}(L\unb\unz_{\e,s}(\unc_{\e,s}(L')))))^{\{-2a-\r\}}\sub
(\io^*(\unc_{\e,s}(L')*\unc_{\e\i,-s}(L)))^{\{0\}}$$
or equivalently, using (a) with $K,K'$ replaced by $b^*\unc_{\e,s}(L')$, $\unc_{\e\i,-s}(L)$:
$$\align&(\io^*(\unc_{1,0}(L\unb\unz_{\e,s}(\unc_{\e,s}(L')))))^{\{-2a-\r\}}\sub
\Hom_{\cc^\boc G_{\e\i,-s}}(\fD(b^*\unc_{\e,s}(L')),\unc_{\e\i,-s}(L))\\&=
\Hom_{\cc^\boc G_{\e,s}}(\fD(b_!\unc_{\e\i,-s}(L)),\unc_{\e,s}(L')).\endalign$$
Using now \cite{\TC, 8.6(c)}, we deduce that we have canonically
$$\Hom_{\cc^\boc\cb^2}(\bold1',L\unb\unz_{\e,s}\unc_{\e,s}L')\sub
\Hom_{\cc^\boc G_{\e,s}}(\fD(b_!\unc_{\e\i,-s}(L),\unc_{\e,s}(L'))$$
where $\bold1'$ is as in \cite{\TC, 8.6} or equivalently (see \cite{\TC, 8.7}):
$$\Hom_{\cc^\boc\cb^2}(\fD(b'_!L),\unz_{\e,s}\unc_{\e,s}L')\sub
\Hom_{\cc^\boc G_{\e,s}}(\fD(b_!\unc_{\e\i,-s}(L),\unc_{\e,s}(L'))$$
where $b':\cb^2@>>>\cb^2$ is $(B,B')\m(B',B)$. We now set ${}^1L=\fD(b'_!L)$ and note that 
$$\fD(b_!\unc_{\e\i,-s}(L))=\fD(\unc_{\e,s}(b'_!L))=\unc_{\e,s}(\fD(b'_!L))=\unc_{\e,s}({}^1L),$$
see 1.3, 1.10(a). We obtain
$$\Hom_{\cc^\boc\cb^2}({}^1L,\unz_{\e,s}\unc_{\e,s}L')\sub
\Hom_{\cc^\boc G_{\e,s}}(\unc_{\e,s}({}^1L),\unc_{\e,s}(L'))\tag b$$
for any ${}^1L,L'\in\cc^\boc_0\cb^2$.

We have the following result which is a weak form of an adjunction formula, of which the full form will be 
proved in 6.6.

\proclaim{Proposition 5.3}Let $s\in\ZZ$. For any ${}^1L,L'\in\cc^\boc_0\cb^2$ we have canonically
$$\Hom_{\cc^\boc\cb^2}({}^1L,\unz_{\e,s}\unc_{\e,s}(L'))
=\Hom_{\cc^\boc G_{\e,s}}(\unc_{\e,s}({}^1L),\unc_{\e,s}(L'))\tag a$$
\endproclaim
We can assume that ${}^1L=\LL_z,L'=\LL_u$ where $z,u\in\boc$. By 4.9(a) and 1.8(b), both sides of the 
inclusion 5.2(b) have dimension $\sum_{y\in\boc}\t(t_{y\i}t_zt_{\e(y)}t_{u\i})$. Hence that inclusion is an 
equality. The proposition is proved. (The case where $\e=1,s=0$ is treated in \cite{\TC, 8.9}.)

\head 6. Equivalence of $\cc^\boc G_{\e,s}$ with the $\e$-centre of $\cc^\boc\cb^2$\endhead
\subhead 6.1\endsubhead
For $\e'\in\fA$ such that $\e'(\boc)=\boc$ and $s,s'\in\ZZ$, the bifunctor 
$\cc^\boc_0G_{\e,s}\T\cc^\boc_0G_{\e',s'}@>>>
\cc^\boc_0G_{\e\e',s+s'}$, $K,K'\m K\unst K'$ in 3.6 defines a bifunctor 
$\cc^\boc G_{\e,s}\T\cc^\boc G_{\e',s'}@>>>\cc^\boc G_{\e\e',s+s'}$ denoted again by $K,K'\m K\unst K'$ as 
follows. Let $K\in\cc^\boc G_{\e,s}$, $K'\in\cc^\boc G_{\e',s'}$; we choose mixed structures of pure weight 
$0$ on $K,K'$ (this is possible if $s_0$ in 0.3 is large enough), we define 
$K\unst K'\in\cc^\boc_0G_{\e\e',s+s'}$ as in 3.6 in terms of 
these mixed structures and we then disregard the mixed structure on $K\unst K'$. The resulting object of 
$\cc^\boc G_{\e\e',s+s'}$ is denoted again by $K\unst K'$; it is independent of the choices made.

In the same way, the bifunctor $\cc^\boc_0\cb^2\T\cc^\boc_0\cb^2@>>>\cc^\boc_0\cb^2$, $L,L'\m L\unb L'$ gives
rise to a bifunctor $\cc^\boc\cb^2\T\cc^\boc\cb^2@>>>\cc^\boc\cb^2$ denoted again by $L,L'\m L\unb L'$; the 
functor $\unc_{\e,s}:\cc^\boc_0\cb^2@>>>\cc^\boc_0G_{\e,s}$ gives rise to a functor 
$\cc^\boc\cb^2@>>>\cc^\boc G_s$ denoted again by $\unc_{\e,s}$ (it is again called {\it truncated 
induction}); the functor $\unz_{\e,s}:\cc^\boc_0G_{\e,s}@>>>\cc^\boc_0\cb^2$ gives rise to a functor 
$\cc^\boc G_{\e,s}@>>>\cc^\boc\cb^2$ denoted again by $\unz_{\e,s}$ (it is again called {\it 
truncated restriction}).

The operation $K\unst K'$ is again called {\it truncated convolution}. It has a canonical associativity 
isomorphism (deduced from that in 3.7) which satisfies the pentagon property. 

The operation $L\unb L'$ makes $\cc^\boc\cb^2$ into a monoidal abelian category (see also \cite{\CEL}) which
has a unit object (see \cite{\TC, 9.2}) and is rigid (see \cite{\TC, 9.3}).

Note that $L\m{}^\e L$ (see 1.6) can be regarded as a functor $\cc^\boc\cb^2@>>>\cc^\boc\cb^2$.

\subhead 6.2\endsubhead
Extending slightly a definition in \cite{\MU, 3.1} we define a {\it $\e$-half braiding} for an object
$\cl\in\cc^\boc\cb^2$ as a collection $e_\cl=\{e_{\cl}(L);L\in\cc^\boc\cb^2\}$ where 
$e_{\cl}(L)$ are isomorphisms ${}^\e L\unb\cl@>\si>>\cl\unb L$ such that (i),(ii) below hold:

(i) If $L@>t>>L'$ is any morphism in $\cc^\boc\cb^2$ then the diagram
$$\CD  
{}^\e L\unb\cl@>e_{\cl}(L)>>\cl\unb L\\
@V{}^et\unb1 VV        @V1\unb t VV\\
{}^\e L'\unb\cl@>e_{\cl}(L')>>\cl\unb L'\endCD$$
is commutative.

(ii) If $L,L'\in\cc^\boc\cb^2$ then $e_{\cl}(L\unb L):{}^\e(L\unb L')\unb\cl@>>>\cl\unb(L\unb L')$ is equal
to the composition
$${}^\e L\unb{}^\e L'\unb\cl@>1\unb e_{\cl}(L')>>{}^\e L\unb\cl\unb L'>e_{\cl}(L)\unb1>>\cl\unb L\unb L'.$$
When $\e=1$, this reduces to the definition of a half-braiding for $\cl$ given in \cite{\MU, 3.1}.
Let $\cz_\e^\boc$ the category whose objects are the pairs consisting of an object $\cl$ of $\cc^\boc\cb^2$
and an $\e$-half braiding for $\cl$. For $(\cl,e_\cl),(\cl',e_{\cl'})$ in $\cz_\e^\boc$ we define
$\Hom_{\cz_\e^\boc}((\cl,e_\cl),(\cl',e_{\cl'}))$ to be the vector space consisting of all
$t\in\Hom_{\cc^\boc\cb^2}(\cl,\cl')$ such that for any $L\in\cc^\boc\cb^2$ the diagram
$$\CD  
{}^\e L\unb\cl@>e_{\cl}(L)>>\cl\unb L\\
@V 1\unb t VV        @Vt\unb1VV\\
{}^\e L\unb\cl'@>e_{\cl'}(L)>>\cl'\unb L\endCD$$
is commutative. We say that $\cz_\e^\boc$ is the $\e$-centre of $\cc^\boc\cb^2$. (When $\e=1$, it reduces to
the centre of $\cc^\boc\cb^2$, see \cite{\MU, 3.2}.)

If $s\in\ZZ$ and $K\in\cc^\boc G_s$ then the isomorphisms 2.9(a) provide an $\e$-half braiding for 
$\unz_{\e,s}(K)\in\cc^\boc\cb^2$ so that $\unz_{\e,s}(K)$ can be naturally viewed as an object of 
$\cz_\e^\boc$ denoted by $\ov{\unz_{\e,s}(K)}$. (Note that 2.9 is stated in the mixed category but, as 
above, it implies the corresponding result in the unmixed category.) Then $K\m\ov{\unz_{\e,s}(K)}$ is a 
functor $\cc^\boc G_{\e,s}@>>>\cz_\e^\boc$. We have the following result.

\proclaim{Theorem 6.3} Let $s\in\ZZ$. The functor $\cc^\boc G_{\e,s}@>>>\cz_\e^\boc$,
$K\m\ov{\unz_{\e,s}(K)}$ is an equivalence of categories. 
\endproclaim
When $\e=1,s=0$ this reduces to \cite{\TC, 9.5}. The general case will be proved in 6.5.

Note that, when combined with 1.6(b), the theorem yields for any $F\in G_{\e,s}$ (with $s>0$) a category 
equivalence
$$\Rep^\boc(G^F)@>\si>>\cz_\e^\boc.\tag a$$

\subhead 6.4\endsubhead
By a variation of a general result on semisimple rigid monoidal categories in \cite{\ENO, Proposition 5.4}, 
for any $L\in\cc^\boc\cb^2$ one can define directly an $\e$-half braiding on the object 
$I_\e(L):=\op_{y\in\boc}\LL_y\unb L\unb\LL_{\e(y)\i}$ of $\cc^\boc\cb^2$ such that, denoting by 
$\ov{I_\e(L)}$ the corresponding object of $\cz_\e^\boc$, we have canonically
$$\Hom_{\cc^\boc\cb^2}(L,\cl)=\Hom_{\cz_\e^\boc}(\ov{I_\e(L)},(\cl,e_{\cl})\tag a$$
for any $(\cl,e_{\cl})\in\cz_\e^\boc$. 

The $\e$-half braiding on $I_\e(L)$ can be described as follows: for any $L'\in\cc^\boc\cb^2$ we have 
canonically
$$\align&{}^\e L'\unb I\e(L)=\op_{y\in\boc}{}^\e L'\unb\LL_y\unb L\unb\LL_{\e(y)\i}=
\op_{y,z\in\boc}\Hom_{\cc^\boc\cb^2}(\LL_z,{}^\e L'\unb\LL_y)\ot\LL_z\unb L\unb\LL_{\e(y)\i}\\&
=\op_{y,z\in\boc}\Hom_{\cc^\boc\cb^2}(\LL_{y\i},\LL_{z\i}\unb{}^\e L')\ot\LL_z\unb L\unb\LL_{\e(y)\i}\\&
=\op_{y,z\in\boc}\Hom_{\cc^\boc\cb^2}(\LL_{\e(y)\i},\LL_{\e(z)\i}\unb L')\ot\LL_z\unb L\unb\LL_{\e(y)\i}\\&
=\op_{z\in\boc}\LL_z\unb L\unb\LL_{\e(z)\i}\unb L'=I_\e(L)\unb L'.\endalign$$
By a variation of results in \cite{\MU, 3.3}, \cite{\ENO, 2.15}, we see that $\cz_\e^\boc$ is a semisimple 
$\bbq$-linear category with finitely many simple objects up to isomorphism. Note that

(b) {\it if $\s=(\cl,e_\cl)$ is a simple object of $\cz_\e^\boc$ then $\s$ is a summand of 
$\ov{I_\e(\LL_z)}$ for some $z\in\boc$.}
\nl
Indeed, let $z\in\boc$ be such that $\LL_z$ is a summand of $\cl$ in $\cc^\boc\cb^2$; then by (a), $\s$ is a
summand of $\ov{I_\e(\LL_z)}$.

\subhead 6.5\endsubhead
Let $s\in\ZZ$. For $x\in\boc$ we have canonically $\unz_{\e,s}\unc_{\e,s}\LL_x=I_\e(\LL_x)$ as objects of 
$\cc^\boc\cb^2$, see Theorem 4.8. This identification is compatible with the $\e$-half braidings (see 6.2, 
6.4). (When $\e=1, s=0$ this follows from the last commutative diagram in \cite{\TC, 7.9}; in the general 
case we have an analogous commutative diagram, which is established using the results in Section 4.) It 
follows that
$$\ov{\unz_{\e,s}\unc_{\e,s}\LL_x}=\ov{I_\e(\LL_x)}.\tag a$$
Using this and 6.4(a) with $\cl=\ov{\unz_{\e,s}\unc_{\e,s}\tL}$, $\tL\in\cc^\boc\cb^2$, we see that
$$\Hom_{\cc^\boc\cb^2}(\LL_x,\unz_{\e,s}\unc_{\e,s}\tL)=\Hom_{\cz_\e^\boc}
(\ov{\unz_{\e,s}\unc_{\e,s}\LL_x},\ov{\unz_{\e,s}\unc_{\e,s}\tL}).$$
Combining this with 5.3 we obtain for $\tL=\LL_{x'}$ (with $x'\in\boc$):
$$\AA_{x,x'}=\AA'_{x,x'}\tag b$$
where
$$\AA_{x,x'}=\Hom_{\cc^\boc G_{\e,s}}(\unc_{\e,s}(\LL_x),\unc_{\e,s}(\LL_{x'})),
\AA'_{x,x'}=\Hom_{\cz_\e^\boc}(\ov{\unz_{\e,s}\unc_{\e,s}\LL_x},\ov{\unz_{\e,s}\unc_{\e,s}\LL_{x'}}).$$
Note that the identification (b) is induced by the functor $K\m\ov{\unz_{\e,s}(K)}$.
Let $\AA=\op_{x,x'\in\boc}\AA_{x,x'},\AA'=\op_{x,x'\in\boc}\AA'_{x,x'}$. Then from (b) we have $\AA=\AA'$. 
Note that this identification is compatible with the obvious algebra structures of $\AA,\AA'$.

For any $A\in CS_{\e,s,\boc}$ we denote by $\AA_A$ the set of all $f\in\AA$ such that for any $x,x'$, the 
$(x,x')$-component of $f$ maps the $A$-isotypic component of $\unc_{\e,s}(\LL_x)$ to the $A$-isotypic 
component of $\unc_{\e,s}(\LL_{x'})$ and any other isotypic component of $\unc_{\e,s}(\LL_x)$ to $0$. Then 
$\AA=\op_{A\in CS_{\e,s,\boc}}\AA_A$ is the decomposition of $\AA$ into a sum of simple algebras (each 
$\AA_A$ is $\ne0$ since, by 1.5(b) and 1.8(a), any $A$ is a summand of some $\unc_{\e,s}(\LL_x)$).

Let $\fS$ be a set of representatives for the isomorphism classes of simple objects of $\cz_\e^\boc$. For 
any $\s\in\fS$ we denote by $\AA'_\s$ the set of all $f'\in\AA'$ such that for any $x,x'$, the 
$(x,x')$-component of $f'$ maps the $\s$-isotypic component of $\ov{\unz_{\e,s}\unc_{\e,s}(\LL_x)}$ to the 
$\s$-isotypic component of $\ov{\unz_{\e,s}\unc_{\e,s}(\LL_{x'})}$ and any other isotypic component of 
$\ov{\unz_{\e,s}\unc_{\e,s}(\LL_x)}$ to $0$. Then $\AA'=\op_{\s\in\fS}\AA'_\s$ is the decomposition of 
$\AA'$ into a sum of simple algebras (each $\AA'_\s$ is $\ne0$ since any $\s$ is a summand of some 
$\ov{\unz_{\e,s}\unc_{\e,s}(\LL_z)}$ (we use 6.4(b), 6.5(a)).

Since $\AA=\AA'$, from the uniqueness of decomposition of a semisimple algebra as a direct sum of simple
algebras, we see that there is a unique bijection $CS_{\e,s,\boc}\lra\fS$, $A\lra\s_A$ such that
$\AA_A=\AA'_{\s_A}$ for any $A\in CS_{\e,s,\boc}$. From the definitions we now see that for any 
$A\in CS_{\e,s,\boc}$ we have $\ov{\unz_{\e,s}A}\cong\s_A$. Therefore Theorem 6.3 holds.

\proclaim{Theorem 6.6} Let $s\in\ZZ$. Let $L\in\cc^\boc\cb^2$, $K\in\cc^\boc G_{\e,s}$. We have canonically
$$\Hom_{\cc^\boc\cb^2}(L,\unz_{\e,s}(K))=\Hom_{\cc^\boc G_{\e,s}}(\unc_{\e,s}(L),K).\tag a$$
Moreover, in $\cc^\boc\cb^2$ we have $\unz_{e,s}(K)\cong\op_{z\in\boc^0}\LL_z^{\op m_z}$ where $\boc^0$ is 
as in 1.5 and $m_z\in\NN$.
\endproclaim
From 6.3, 6.5, we see that 
$$\align&\Hom_{\cc^\boc G_{\e,s}}(\unc_{\e,s}(L),K)=\Hom_{\cz_\e^\boc}(\ov{\unz_{\e,s}\unc_{\e,s}(L)},
\ov{\unz_{\e,s}K})\\&=\Hom_{\cz_\e^\boc}(\ov{I_\e(L)},\ov{\unz_{\e,s}K}).\endalign$$
Using 6.4(a) we see that 
$$\Hom_{\cz_\e^\boc}(\ov{I_\e(L)},\ov{\unz_{\e,s}K})=\Hom_{\cc^\boc\cb^2}(L,\unz_{\e,s}(K))$$ 
and (a) follows. To prove the second assertion of the theorem it is enough to show that for any 
$z\in\boc-\boc^0$ we have $\Hom_{\cc^\boc\cb^2}(\LL_z,\unz_{\e,s}(K))=0$; by (a), it is enough to show that 
$\unc_{\e,s}(\LL_z)=0$ and this follows from 1.5(c). (The case where $\e=1,s=0$ is just \cite{\TC, 9.8}.)

\subhead 6.7\endsubhead
Let $s\in\ZZ$. For $K\in\cc^\boc G_{\e,s}$ we have canonically
$$\fD(\unz_{\e,s}(\fD(K)))=\unz_{\e,s}(K).\tag a$$
When $\e=1,s=0$ this is proved in \cite{\TC, 9.9}. The proof in the general case is entirely similar; it uses
6.6(a), 1.10(a).

\subhead 6.8\endsubhead
In this subsection we assume that $\e=1$. The monoidal structure on $\cc^\boc\cb^2$ induces a monoidal 
structure on $\cz_1^\boc$. Moreover, the category
$$\sqc_{s\in\ZZ}\cc^\boc G_{1,s}=\cc^\boc G_{1,0}\sqc_{s\in\ZZ;s\ne0}\Rep^\boc(G^{F_0^s})\tag a$$
(see 1.6(b)) has a monoidal structure given by truncated convolution, see 6.1. Moreover, 6.3 provides a
functor from (a) to $\cz_1^\boc$ which is an equivalence when restricted to any $\cc^\boc G_s$. This functor
is compatible with the monoidal structures (this can be deduced from 3.8 and from the fact that the 
monoidal structure of $\cz_1^\boc$ is equivalent to its opposite). Note that $\cc^\boc G_{1,0}$ is a 
monoidal subcategory of (a), whose unit object, described in \cite{\TC, 9.10}, is also a unit object for
the monoidal category (a).

\subhead 6.9\endsubhead
The functor $L\m{}^\e L$ from $\cc^\boc\cb^2$ into itself induces a functor $\cz_\e^\boc@>>>\cz_\e^\boc$
which carries any simple object $(L,e_L)$ of $\cz_\e^\boc$ into an object isomorphic to $(L,e_L)$; this 
follows from 2.10(c), using Theorem 6.3.

\subhead 6.10\endsubhead
Let $s\in\ZZ$.
For any $A\in CS_{\e,s,\boc}$ and any $x\in\boc$ we denote by $n_{A,x}$ the multiplicity of $A$ in 
$\c_{\e,s}\LL_x\in\cc^\boc G_{\e,s}$. From Theorem 6.3 and its proof we see that if $\s$ is the simple 
object of $\cz_\e^\boc$ corresponding to $A$, then $n_{A,x}$ is equal to the multiplicity of $\s$ in 
$\ov{I_\e(\LL_x)}\in\cz_\e^\boc$. In particular, $n_{A,x}$ is independent of $s$.

\head 7. Relation with Soergel bimodules\endhead
\subhead 7.1\endsubhead
Let $R$ be the algebra of polynomials functions on a fixed reflection representation of $W$ (over $\bbq$).
Then for each $x\in W$, the indecomposable Soergel graded $R$-bimodule $B_x$ is defined as in 
\cite{\SO, 6.16}. Let $C_{\boc}$ be the category of graded $R$-bimodules wich are isomorphic to finite 
direct sums of graded $R$-bimodules of the form $B_x$ ($x\in\boc$) without shift. There is a well defined
functor $M\m{}^\e M$ from $C_\boc$ to $C_\boc$ which is linear and satisfies ${}^\e B_x=B_{\e\i(x)}$ for
$x\in\boc$. Now $C_\boc$ has a natural
monoidal structure (see \cite{\TC, 10.1} defined purely in terms of $R,W,\boc$. (Its definition makes use
of the results in \cite{\EW}.) From the definition we see that $C_{\boc}$ is equivalent to $\cc^\boc\cb^2$ as
monoidal categories so that $M\m{}^\e M$ corresponds to $L\m{}^\e L$ from $\cc^\boc\cb^2$ to itself.
Then the $\e$-centre of $C_\boc$ is defined as in 6.2. It is naturally equivalent to $\cz_\e^\boc$.
Thus we can restate Theorem 6.3 as follows.

(a) {\it For any $s\in\ZZ$, the category $\cc^\boc G_{\e,s}@>>>\cz_\e^\boc$ is naturally equivalent to
the $\e$-centre of the monoidal category $C_\boc$.}
\nl
This, combined with 1.6(b), shows that for $F\in G_{\e,s}$ (with $s>0$), the category $\Rep^\boc(G^F)$ is
equivalent to the $\e$-centre of the monoidal category $C_\boc$; thus, the set of simple objects of 
$\Rep^\boc(G^F)$ is not only independent of $s$ but also independent of the characteristic of $\kk$, since 
the $\e$-centre of $C_\boc$ is so. (Here we identify $\bbq$ with the complex numbers.)

\subhead 7.2\endsubhead
As mentioned in \cite{\TC, 10.1}, the definition of the monoidal category $C_\boc$ makes sense even when $W$
is replaced by any (say finite, irreducible) Coxeter group and $\boc$ is a two-sided cell in $W$.
Assume now that $\e:W@>>>W$ is an automorphism of $W$ which leaves stable the set of 
simple reflections and leaves stable $\boc$. Then the definition of the $\e$-centre of $C_\boc$ makes sense
even if $W$ is noncrystallographic.
We expect that the indecomposable objects of the $\e$-centre of $C_\boc$ are in bijection with the 
``unipotent characters'' associated to $W,\e,\boc$ in \cite{\COX}. (For $\e=1$ this expectation has already 
been stated in \cite{\TC, 10.1}.)


\widestnumber\key{ENO}
\Refs
\ref\key\BFO\by R.Bezrukavnikov, M.Finkelberg and V.Ostrik\paper Character D-modules via Drinfeld center of
Harish-Chandra bimodules\jour Invent. Math.\vol188\yr2012\pages589-620\endref
\ref\key\DL\by P.Deligne and G.Lusztig\paper Representations of reductive groups over finite fields\jour
Ann. Math.\vol103\yr1976\pages103-161\endref
\ref\key\EW\by B.Elias and G.Williamson\paper The Hodge theory of Soergel bimodules\jour arxiv:1212.0791
\endref
\ref\key\ENO\by P.Etingof, D.Nikshych and V.Ostrik\paper On fusion categories\jour Ann. Math.\vol162\yr2005
\pages581-642\endref
\ref\key\ORA\by G.Lusztig\book Characters of reductive groups over a finite field
\bookinfo Ann. Math. Studies 107\publ Princeton U.Press\yr1984\endref 
\ref\key\CSII\by G.Lusztig\paper Character sheaves,II\jour Adv. Math.\vol57\yr1985\pages226-265\endref
\ref\key\CSIII\by G.Lusztig\paper Character sheaves,III\jour Adv. Math.\vol57\yr1985\pages266-315\endref
\ref\key\CSIV\by G.Lusztig\paper Character sheaves,IV\jour Adv. Math.\vol59\yr1986\pages1-63\endref
\ref\key\COX\by G.Lusztig\paper Coxeter groups and unipotent representations\jour Ast\'erisque\vol212\yr1993
\pages191-203\endref
\ref\key\CDGVII\by G.Lusztig\paper Character sheaves on disconnected groups, VII\jour Represent.Th.\vol9
\yr2005\pages209-266\endref
\ref\key\CDGIX\by G.Lusztig\paper Character sheaves on disconnected groups, IX\jour Represent.Th.\vol10
\yr2006\pages353-379\endref
\ref\key\CDGX\by G.Lusztig\paper Character sheaves on disconnected groups, X\jour Represent.Th.\vol13
\yr2009\pages82-140\endref
\ref\key\CEL\by G.Lusztig\paper Cells in affine Weyl groups and tensor categories\jour Adv. Math.\vol129\yr
1997\pages85-98\endref
\ref\key\HEC\by G.Lusztig\book Hecke algebras with unequal parameters\bookinfo CRM Monograph Ser.18\publ
Amer. Math. Soc.\yr2003\endref
\ref\key\RES\by G.Lusztig\paper Restriction of a character sheaf to conjugacy classes\jour arxiv:1204.3521
\endref
\ref\key\TC\by G.Lusztig\paper Truncated convolution of character sheaves\jour arxiv:1308.1082\endref
\ref\key\MU\by M.M\"uger\paper From subfactors to categories and topology II. The quantum double of tensor
categories and subfactors\jour J. Pure Appl. Alg.\vol180\yr2003\page159-219\endref
\ref\key\SO\by W.Soergel\paper Kazhdan-Lusztig Polynome und unzerlegbare Bimoduln \"uber Polynomringen\jour 
J. Inst. Math. Jussieu\vol6\yr2007\pages501-525\endref

\endRefs
\enddocument